\documentclass[<option>]{elsarticle}
\usepackage{tipa}
\usepackage{textcomp}

\usepackage{amsfonts,mathrsfs}
\usepackage{latexsym,amssymb}
\usepackage{amsmath, amsbsy}
\usepackage{amsopn, amstext}
\usepackage{graphicx, color, epstopdf}
\usepackage{threeparttable}
\usepackage[colorlinks]{hyperref}
\usepackage{array}
\usepackage{bm}
\usepackage{multirow}
\usepackage{booktabs}
\usepackage{diagbox}
\usepackage{subcaption}
\usepackage{graphics}
\usepackage{indentfirst}
\usepackage{enumitem}
\usepackage{multirow}
\usepackage{float}
\usepackage{tikz}
\usepackage{overpic}
\usetikzlibrary{graphs}
\usepackage{tabularx}
\usepackage{booktabs,arydshln} 
\usepackage{makecell} 
\usepackage{amsmath}
\allowdisplaybreaks[0]

\usepackage{algorithm}
\usepackage[noend]{algpseudocode}

\usepackage{pgflibraryarrows}
\usepackage{pgflibrarysnakes}
\usepackage{hyperref}
\hypersetup{hypertex=true,
colorlinks=true,
linkcolor=blue,
anchorcolor=blue,
citecolor=blue}

\def\R{\mathbb{R}}

\numberwithin{equation}{section}
\numberwithin{figure}{section}
\numberwithin{table}{section}

\allowdisplaybreaks[4]

\newtheorem{remark}{Remark}[section]

\numberwithin{equation}{section}

\def\R{\mathbb{R}}

\def\N{\mathcal{N}}

\def\I{\mathcal{I}}
\def\L{\mathcal{L}}

\def\V{\mathcal{V}}
\def\U{\mathcal{U}}

\def\G{\mathcal{G}}
\def\|{\Vert}
\newcommand{\vx}{\bm{x}}

\def\ft{f^{tar}}
\def\varphit{\varphi^{tar}}
\def\gt{g^{tar}}

\newcommand{\beq}{\begin{eqnarray}}
\newcommand{\eeq}{\end{eqnarray}}
\newcommand{\beqs}{\begin{eqnarray*}}
\newcommand{\eeqs}{\end{eqnarray*}}

\begin{document}
\begin{frontmatter}
\title{\bf Data Generation-based Operator Learning for Solving Partial Differential Equations on Unbounded Domains}

\author[1]{Jihong Wang}
\ead{jhwang@zhejianglab.com}

\author[2]{Xin Wang}
\ead{xinwang2021@whu.edu.cn}

\author[3]{Jing Li \corref{cor1}}
\ead{lijing@zhejianglab.com}

\author[1]{Bin Liu \corref{cor1}}
\ead{liubin@zhejianglab.com}
\cortext[cor1]{Corresponding author.}

\address[1]{Research Center for Frontier Fundamental Studies, Zhejiang Lab, Hangzhou 311121, China}

\address[2]{School of Mathematics and Statistics, Wuhan University, Wuhan 430072, China}

\address[3]{Research Center for Intelligent Equipment, Zhejiang Lab, Hangzhou 311121, China}

\begin{abstract}

  Wave propagation problems are typically formulated as partial differential equations (PDEs) on unbounded domains to be solved.
  The classical approach to solving such problems involves truncating them to problems on bounded domains by designing the artificial boundary conditions or perfectly matched layers, which typically require significant effort, and the presence of nonlinearity in the equation makes such designs even more challenging. 
  Emerging deep learning-based methods for solving PDEs, with the physics-informed neural networks (PINNs) method as a representative, still face significant challenges when directly used to solve PDEs on unbounded domains. Calculations performed in a bounded domain of interest without imposing boundary constraints can lead to a lack of unique solutions thus causing the failure of PINNs.
  In light of this, this paper proposes a novel and effective data generation-based operator learning method for solving PDEs on unbounded domains. 
  The key idea behind this method is to generate high-quality training data.
  Specifically, we construct a family of approximate analytical solutions to the target PDE based on its initial condition and source term.
  Then, using these constructed data comprising exact solutions, initial conditions, and source terms, we train an operator learning model called MIONet, which is capable of handling multiple inputs, to learn the mapping from the initial condition and source term to the PDE solution on a bounded domain of interest. 
  Finally, we utilize the generalization ability of this model to predict the solution of the target PDE.
  The effectiveness of this method is exemplified by solving the wave equation and the Schr\"odinger equation defined on unbounded domains. 
  More importantly, the proposed method can deal with nonlinear problems, which has been demonstrated by solving Burgers' equation and Korteweg-de Vries (KdV) equation.  The code is available at \url{https://github.com/ZJLAB-AMMI/DGOL}.

\end{abstract}
  
\begin{keyword}
  Scientific machine learning \sep Operator learning\sep  Unbounded domain\sep Nonlinear PDEs
  \end{keyword}

 \end{frontmatter}

  \section{Introduction}
  Real-world wave propagation problems in various fields, such as acoustics, aerodynamics, solid geophysics, oceanography, meteorology, and electromagnetics, are commonly described by partial differential equations (PDEs) on unbounded (or very large) domains. 
  However, solving these PDEs numerically poses challenges due to the infinite domains involved.
  Standard domain-based numerical methods like finite difference and finite element methods are not directly applicable to solving such unbounded PDEs. The main reason is that these methods rely on discretizing the domain into a finite set of points or elements, leading to a reduced algebraic system with a finite number of degrees of freedom. However, when dealing with unbounded domains, these approaches lead to an algebraic system with an infinite number of degrees of freedom that cannot be effectively solved.
  
  Two popular techniques to deal with the problems on unbounded domains are the artificial boundary method (ABM) \cite{han2013artificial,antoine2008review,givoli2004high} and perfectly matched layer (PML) \cite{berenger1994perfectly,pled2022review}.
  The ABM involves designing suitable absorbing/artificial boundary conditions (ABCs) that are satisfied by the solution of the original problem on the artificial boundaries. This approach reduces the original unbounded problem to a well-posed boundary value problem on the bounded computational domains of interest. The key ingredient of ABM is the construction of the ABCs. Based on the Fourier series expansion, Laplace transform, $z$-transform, Pad\'e approximation, continued fraction expansion, and other techniques, exact or approximated ABCs are designed for various linear PDEs. However, it remains challenging to construct suitable ABCs for many nonlinear equations, even for some fundamental equations like the nonlinear Schr\"odinger equation \cite{han2013artificial, zheng2006exact}. 
  The PML method is used as absorbing layers that effectively eliminate reflections for all incident waves, regardless of their frequency and angle. PMLs have gained widespread usage due to their computational efficiency, ease of implementation, applicability to complex geometries, and high absorption accuracy. However, a major drawback of the PML method is its numerical instabilities in time-domain simulations for some wave propagation problems in anisotropic and/or dispersive media \cite{becache2003stability,becache2018analysis,duru2012well,loh2009fundamental}. Several methods have been proposed to avoid and remove growing waves and to improve the stability and accuracy of PML formulations. Nevertheless, for some specific situations, such as linear elastodynamic equations in arbitrary anisotropic elastic media, developing a stable PML formulation remains an open problem \cite{pled2022review}.
  
\begin{figure}[ht]
    \centering
    \includegraphics[width=0.4\linewidth]{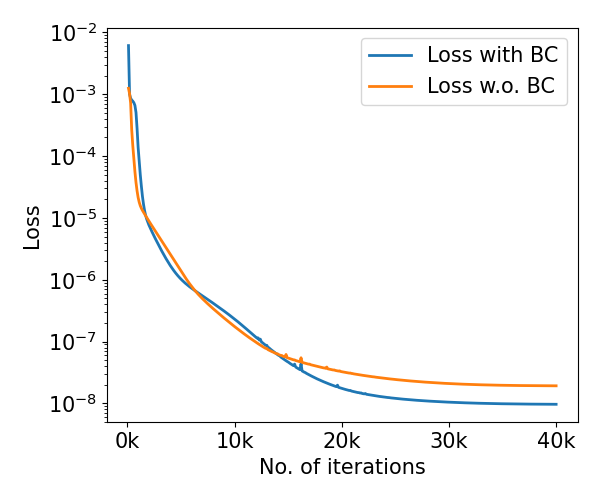}
    \includegraphics[width=0.4\linewidth]{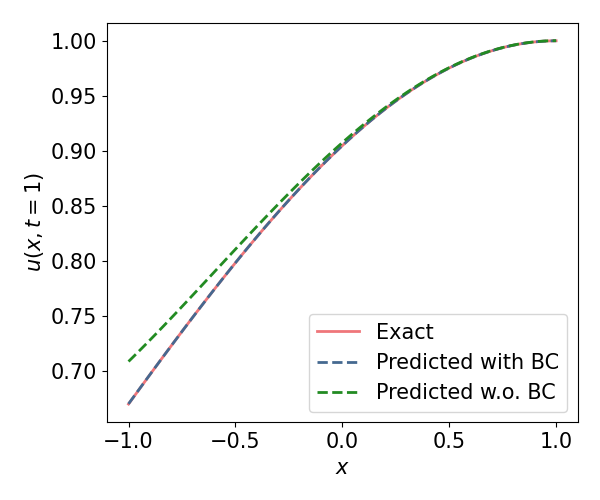}
    \caption{\textit{Solving a 1D Burgers' equation using PINNs with and without boundary conditions:} Left: Training loss. Right: Exact solution and predicted solutions at time $t=1$. 
    It is evident that in the absence of boundary conditions, the PINN solution exhibits a substantial error even after achieving training stability. (The experimental details can be found in the Appendix D.)}
    \label{PINN}
  \end{figure}

  Recently, the rapid progress in deep learning has driven the development of solution techniques for PDEs. Learning-based PDE approaches can fall into two categories in terms of the objects approximated by neural networks (NN), i.e., the solution and the solution mapping. 
  The typical methods of the first category include physics-informed neural networks (PINNs) \cite{karniadakis2021physics}, Deep Galerkin Method (DGM) \cite{sirignano2018dgm}, Deep Ritz Method \cite{yu2018deep}, Weak Adversarial Network (WAN) \cite{zang2020weak}, etc. These methods train neural networks by sampling points within bounded computational domains and minimizing the physical loss at these sampled points. However, when dealing with PDEs on unbounded domains, they become impractical as it would require sampling training data points across the entire unbounded domains. On the other hand, if data points are only sampled within a bounded domain without imposing boundary conditions, the problem will become ill-posed. This can result in the neural network training process failing to converge or converging to an incorrect solution, as illustrated in Figure \ref{PINN}. 

  Methods of the second category use neural networks to learn the solution mapping between infinite-dimensional function spaces (such as the mapping from the initial function space to the solution space). The typical methods include deep operator network (DeepONet) \cite{lu2021learning}, Fourier neural operator (FNO) \cite{li2021fourier}, PDE-Nets \cite{long2018pde,long2019pde}, etc. More related work can be found in \cite{khoo2021solving,zafar2022frame,zhang2022mod,lu2022comprehensive}. A significant advantage of this approach is that once the neural network is trained, it can quickly solve a large number of in-distribution tasks. However, this approach requires a large amount of data to effectively train neural networks, which is costly and challenging in many cases. For the task of solving PDEs on bounded domains, it is generally necessary to use traditional numerical methods to extensively solve PDEs to obtain data for training. For PDEs defined on unbounded domains, the cost of data acquisition is higher because of the need for complicated techniques, such as ABM or PML. Further, for many complex nonlinear equations, there is even no effective numerical method to solve them, resulting in the lack of training data.
  An alternative method, physics-informed DeepONet (PI-DeepONet) \cite{wang2021learning}, is proposed to train the neural networks only relying on the equations and given initial/boundary conditions. As the combination of PINNs and DeepONet, PI-DeepONet adopts the architecture of DeepONet while taking the residual of the equation as the loss function during training just like PINNs. Unfortunately, this strategy does not work for solving PDEs on unbounded domains due to the same reason as PINNs. 
  
  At present, most deep learning-based solvers focus on problems defined on bounded domains, while there is few work addressing the complexities associated with solving PDEs on unbounded domains \cite{feng2022solving,wilson2022new,gao2023failure, gao2023failure1,xia2023spectrally,lin2023binet,lin2023bi,xie2021machine}.
  Xia et al. \cite{xia2023spectrally} blend the adaptive spectral method and PINNs to optimize the traditional numerical spectral schemes for solving PDEs in unbounded domains, while an assumption on the asymptotic spatial behavior is needed. 
  Lin et al. \cite{lin2023binet,lin2023bi} propose BINet, a method that combines boundary integral equations with neural networks to solve PDEs with known fundamental solutions. By leveraging potential theory to transform the original problem into boundary integral equations, this approach can handle problems on both bounded and unbounded domains.
  Gao et al. \cite{gao2023failure, gao2023failure1} propose a failure-informed enrichment adaptive sampling PINN, which is appliable to PDEs on unbounded domains with local behavioral solutions.
  Xie et al. \cite{xie2021machine} are the first to explore the application of DeepONet for solving the PDEs on unbounded domains. They train the neural network in a data-driven setting, where data is obtained by exact solutions or traditional numerical methods. The trained model can be used to predict the solution of interpolation problems. 
  This direct application of the DeepONet is confronted with the issues of high data acquisition costs and ineffectiveness for many challenging problems that traditional numerical methods struggle to solve, as mentioned earlier.
  Therefore, efficiently solving problems on unbounded domains, especially with the assistance of powerful deep learning techniques, is an area that merits further exploration.

  This work aims to develop an effective operator learning-based approach to solve PDEs (including nonlinear cases) defined on unbounded domains. 
  The basic idea of our approach is to generate training data at a low cost and leverage the generalization ability of DeepONet to predict the solution of the PDE on the domain of interest. To ensure that the PDE to be solved falls within the interpolation range of the training data, we generate paired input-output training data that closely approximates the target PDE's initial value and source term.
  The main operations of the approach are as follows: firstly, construct a family of analytical solutions that satisfy the target equation needed to be solved, with initial value and source term designed to closely `approximate' those of the target equation. Then, utilize MIONet \cite{jin2022mionet}, an extension to DeepONet for learning multiple-input operators, to learn the mapping from initial value and source term to the PDE solution. Finally, utilize the learned mapping to directly predict the solution of the target equation on a bounded domain of interest. 
  We test the effectiveness of our method on extensive equations including the second-order wave equation, Burgers' equation, Korteweg-de Vries (KdV) equation, and Schr\"odinger equation. 
  This approach offers several advantages. Firstly, it involves low computational costs for generating training data to train operator models. Secondly, it demonstrates effectiveness in solving nonlinear PDEs on unbounded domains, which are challenging for traditional numerical methods. Lastly, this method is generally applicable to various types of PDEs verified by numerical experiments and only needs to generate corresponding data that meets specific requirements. This flexible approach eliminates the need for designing distinct techniques for different equations, as often required by traditional numerical methods.
  
  The paper is structured as follows. Section 2 provides a brief introduction to the architecture of DeepONet and MIONet. In Section 3, we present a general workflow of our proposed method for solving PDEs on unbounded domains. Next, in Section 4, we demonstrate the effectiveness of our method through extensive numerical examples including both linear and nonlinear cases. Finally, Section 5 concludes the paper with a discussion of our main findings, potential limitations of the proposed method, and future research directions stemming from this study.

  \section{Preliminaries}

  \begin{figure}[H]
      \centering
  \includegraphics[width=1.0\linewidth]{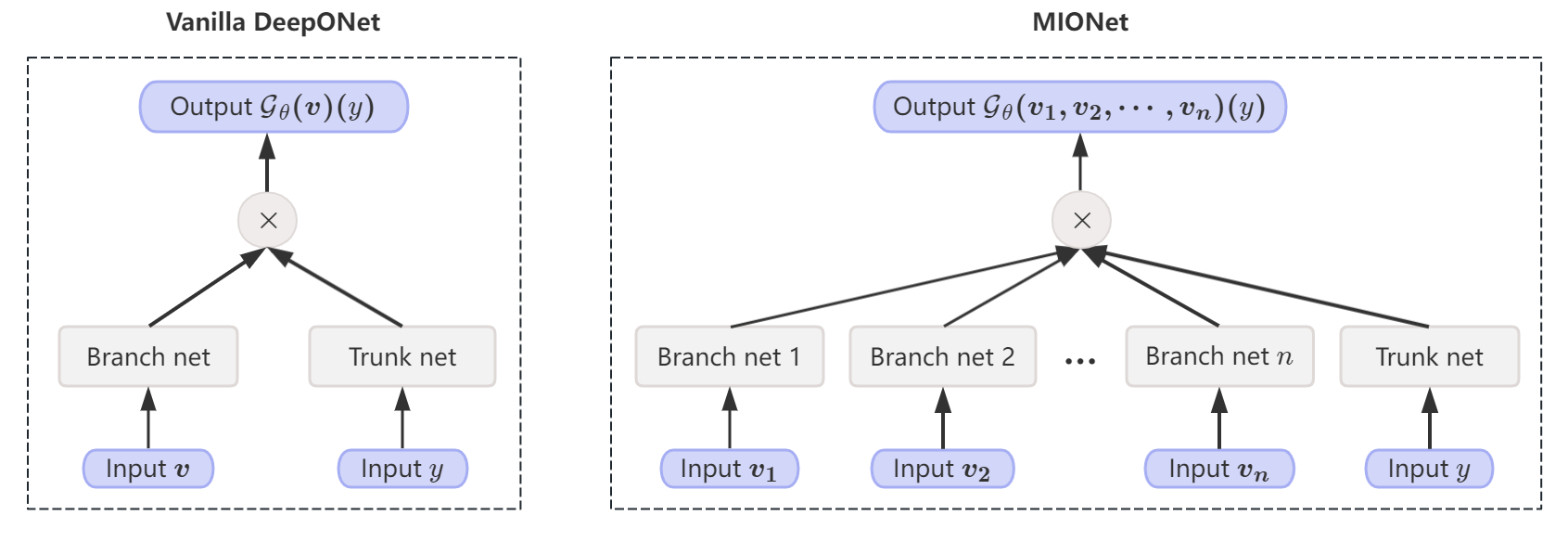}
  \caption{Architectures of vanilla DeepONet and MIONet.}
  \label{NN}
  \end{figure}

  \subsection{DeepONet}
  Based on the universal approximation theorem of Chen $\&$ Chen \cite{chen1995universal}, DeepONet \cite{lu2021learning} is proposed to learn a nonlinear operator
  \beqs
  \G: \V\rightarrow \U, \quad v \mapsto u,
  \eeqs
  where $\V$ and $\U$ are two infinite-dimensional Banach spaces of functions, $v\in \V$ and $u=\G(v)$. The DeepONet learns $\G$ via a neural network denoted by $\G_{\theta}$, where $\theta$ represents all trainable parameters of the network.
  Take $v\in \V$ evaluated at a collection of fixed locations $\{y_j\}_{j=1}^{m}$, i.e., $\bm{v}=[v(y_1),v(y_2),\cdots,v(y_{m})]$ and another family of locations $y=\{y_j\}_{j=1}^{P}$ as inputs, the output of DeepONet is calculated as follows
  \beq
  \G_{\theta}(\bm{v})(y) = S\left( B(\bm{v}) \odot T(y) \right) + b,
  \eeq
  where $B$ is called branch net, and $T$ called trunk net, $\odot$ is the Hadamard product (i.e., element-wise product), $S$ is the summation of all components of a vector, and $b$ is a trainable bias. The architecture of DeepONet is shown in the left plane of Figure \ref{NN} and its detailed dataset structure can be referred to \cite[Remark 2.1]{wang2021learning}.

  \subsection{MIONet}
  The vanilla DeepONet is defined for input functions on a single Banach space. MIONet \cite{jin2022mionet} extends DeepONet to multiple input Banach spaces. Specifically,
  MIONet aims to learn a continuous operator
  \beqs
  \G: \V_1\times \V_2\times \cdots \times \V_n \rightarrow \U, \quad (v_1, v_2,\cdots, v_n) \mapsto u,
  \eeqs
  where $\V_1, \V_2, \cdots, \V_n$ are $n$ different input Banach spaces, $\U$ is the output Banach space, and $v_i\in \V_i$, $u=\G(v_1, \cdots, v_2)$. Similar to the vanilla DeepONet, MIONet approximates $\G$ by
  \beq
  \G_{\theta}(\bm{v_1}, \bm{v_2}, \cdots, \bm{v_n})(y) = S\left( B_1(\bm{v_1}) \odot \cdots \odot B_n(\bm{v_n}) \odot T(y) \right) + b,
  \eeq
  where $\{B_i\}_{i=1}^n$ are $n$ different branch nets, each $\bm{v_i}=[v_i(y^i_1),v_i(y^i_2),\cdots,v_i(y^i_{m})]$ can be calculated at different locations $\{\bm{y^i}\}$.
  The architecture of MIONet is shown at the right plane of Figure \ref{NN}.

  \section{Proposed method}
  In this section, we present the proposed data generation-based operator learning method for solving PDEs on unbounded domains.
  Consider the following time-dependent equation
  \beq
  \L u(x,t) = f(x,t), \quad x\in\R^d,~t>0, \label{Eq}\\
  \I_i u(x,0) = \varphi_i(x), \quad x\in \R^d,~i\in{I}, \label{EqI}
  \eeq
  where $\L$ and $\I_i$ are partial differential operators, and $f(x,t)$ is the source term, $\varphi_i(x)$ is the initial function, $I=\{0\}$ if Eq.\eqref{Eq} is the first-order time evolution equation and $I=\{0,1\}$ when Eq.\eqref{Eq} is the second-order time evolution equation. There are no extra constraints on the initial value and source term. For convenience of exposition, let us refer to the equation that needs to be solved as the target equation, and the corresponding initial value and source term will be referred to as target functions, denoted by $\varphi^{tar}$ and $f^{tar}$, respectively.

\begin{figure}[H]
    \centering
\includegraphics[width=1.0\linewidth]{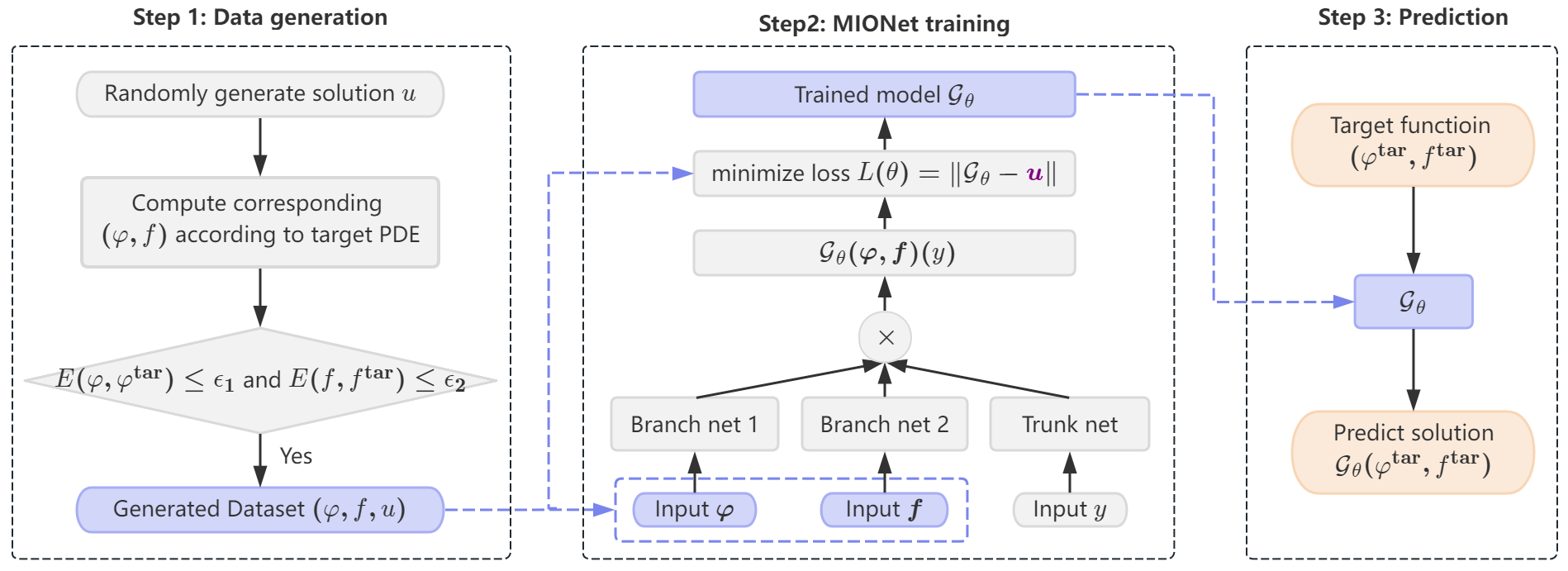}
\caption{Flowchart of the proposed method.}
\label{flow}
\end{figure}

  The proposed method of solving Eq.\eqref{Eq}-\eqref{EqI} can be summarized into the following three steps (see also Figure \ref{flow} for the entire workflow).

  {\bf Step 1: Data generation.}
  
  To generate data for learning a mapping from the initial value and source term to the solution, we construct a family of analytical solutions that satisfy the target governing equation, for which corresponding initial values and source terms approximate the target initial value and the target source term, respectively. We describe how to construct analytical solutions and select training data in detail as follows.

  The core idea of constructing analytical solutions is leveraging the information of known target functions. For example, if the target initial value is zero, we can construct an analytical solution that inherently satisfies this condition, simplifying the model's task to learn the mapping solely from the source term to the solution. Similarly, if the target function is compactly supported, we can generate an analytical solution using the compactly supported basis function accordingly. 
  On the basis of utilizing the information of target functions, the choice of analytical solution in form is still diverse. For instance, we can choose the following solution in the form of variable separation
  \beq\label{u}
  u(x,t) = \sum_{k=0}^K \alpha_k(t) \phi_k(w x+b),
  \eeq
  where $\phi_k(x)$ is some basis function on $x$, such as Hermite function, $\alpha_k(t)$ is basis function on $t$, such as Fourier function. Parameters $K, w,b$ and parameters contained in $\alpha_k(t)$ are randomly generated according to certain rules, such as normal random distributions. 
  After obtaining an analytical solution, substituting it into Eq.\eqref{Eq}-\eqref{EqI}, one obtains the corresponding initial value $\varphi$ and source term $f$. 
  
  To assess the quality of the generated data $\varphi, f$ in comparison to the target function $\varphit, \ft$, we introduce the following evaluation metrics:
  \begin{itemize}
  \item Metric-1 (relative $L^2$ error):
  \beq \label{metric1}
   E_{1}(g,\gt) :=\frac{\|\gt-g\|_2}{\|\gt\|_2},
  \eeq
  where $\|\cdot\|_2$ represents the discrete $L^2$ norm.
  \item Metric-2 (mean square error):
  \beq \label{metric2}
      E_{2}(g,\gt) := \left(\frac{1}{P}\sum_{j=1}^P |\gt(y_j)-g(y_j)|^2\right)^{\frac12},
  \eeq
  where $\{y_j\}_{j=1}^P$ are sample points in a fixed spatial-temporal domain related to the computational domain of interest.
  \end{itemize}
  Generally, we utilize metric-1 to evaluate the generated data. However, when $\gt=0$, metric-1 becomes ineffective, and in such cases, we rely on metric-2 for evaluation. 
  If $E(g,\gt) \le \epsilon$, then the corresponding data is selected to train the model, where $\epsilon$ is a parameter selected based on experience. 

  We remark that the evaluation should better be conducted in a relatively large domain compared to the domain of interest.
  By ensuring that the initial values and source terms closely match the target functions across a broader domain, one can achieve higher-quality generated data.  We investigate the impact of the size of the spatial domain where the data is generated on the final results, as detailed in Appendix B.

  {\bf Step 2: MIONet training.}
  
  MIONet is utilized to learn the mapping from the initial value and source term to the solution of the PDE, i.e., $\mathcal{G_{\theta}}: (\varphi, f) \rightarrow u $.
  We train the network parameters by minimizing the following loss function:
\beq \label{loss_total}
L(\theta) = w_dL_{data}(\theta) + w_pL_{phys}(\theta),
\eeq
where $L_{data}$ is the loss of data
\beq \label{loss_data}
L_{data}(\theta) = \frac{1}{NP} \sum_{n=1}^N \sum_{j=1}^P|\G_{\theta}(\bm{\varphi}^{(n)},\bm{f}^{(n)})(y_j)-u^{(n)}(y_j)|^2, 
\eeq
$L_{phys}$ is the loss of PDE residuals
\beq\label{loss_phys}
L_{phys}(\theta) = \frac{1}{NQ} \sum_{n=1}^N \sum_{j=1}^Q|\L(\G_{\theta}(\bm{\varphi}^{(n)},\bm{f}^{(n)}))(y_j)-f^{(n)}(y_j)|^2,
\eeq
and $w_d$ and $w_p$ are the weights. $\bm{\varphi}^{(n)} = [\varphi^{(n)}(y_1),\varphi^{(n)}(y_2), \cdots, \varphi^{(n)}(y_{m_1})]$ represents the initial function evaluated at a collection of fixed locations $\{y_j\}_{j=1}^{m_1}$, $\bm{f}^{(n)} = [f^{(n)}(y_1),f^{(n)}(y_2), \cdots, f^{(n)}(y_{m_2})]$ represents the source term evaluated at another set of fixed locations $\{y_j\}_{j=1}^{m_2}$, and $u^{(n)}(y_j)$ denotes the solution at points $y_j$ of the problem \eqref{Eq}-\eqref{EqI} with $\varphi^{(n)}$ and $f^{(n)}$. $N$ is the number of input function pairs, $P$ and $Q$ are the number of sampling points on the computational domain used to compute the data loss and the physics loss, respectively. 

In the numerical experiments, we compare the predicted results obtained by training the model using only data loss $L_{data}(\theta)$ and total loss $L(\theta)$. 
Our findings reveal that physics-informed training does not lead to a significant improvement in accuracy. Moreover, it considerably slows down the training process, as detailed in \cite[Appendix C]{wang2021learning}. Consequently, unless explicitly specified otherwise in the experiment, we opted to exclusively employ data loss for training the model.
  
  {\bf Step 3: Prediction.}
  
  Apply the learned mapping operator $\mathcal{G_{\theta}}$ to target functions $\varphit$ and $\ft$ to obtain the solution of the target PDE. 
   
  A natural question arises: how to assess the reliability of predicted solutions? One intuitive approach is to substitute the predicted solution into the equation and compute the residual. A small residual indicates that the predicted solution effectively satisfies the equation, making it highly probable that it is a solution to the original equation.
  However, it is crucial to recognize that when the model is trained with only data, the predicted solution generated by the model may closely resemble the true solution, but the corresponding residual value of its equation is likely to be substantial. Conversely, the physics-informed training typically results in a reduction of the PDE residual for the predicted solution. This observation is exemplified in the numerical example 4.4. Therefore, to ensure more dependable predicted solutions, opting for physics-informed training is a viable consideration, albeit at the expense of increased training time.

  \begin{remark}
  The success of this approach relies on the well-posedness of the equation, particularly its stability concerning the initial value and source term, 
  i.e., there exists a constant $C$, s.t.,
  \beqs
  \|u-u_{\epsilon}\| \le C\|f-f_{\epsilon}\|,
  \eeqs
  where $f_{\epsilon}$ represents the initial function or source function $f$ after a small perturbation, and $u_{\epsilon}$ is the corresponding solution.
  \end{remark}
  \begin{remark}
    While we train the model using analytical solution data, it is important to highlight that, owing to the model's generalization ability, it can predict solutions for problems without analytical solutions. The subsequent numerical examples serve to exemplify this characteristic. 
  \end{remark}

  \section{Numerical results}
  In this section, we evaluate the effectiveness of the proposed method across various types of PDEs, including 1D and 2D wave equations, Burgers' equation, KdV equation and Schr\"odinger equation. All experiments are conducted on an NVIDIA A100 GPU using PyTorch (v1.12.1) framework. The activation function employed in all networks is set to Tanh, and the batch size is fixed at 8192.
  The networks are trained with the Adam optimizer [12] with default settings. 
  Unless specified otherwise, by default, only data loss is used during model training.
  The detailed hyper-parameters for all examples are listed in Appendix A.
  In our evaluation, we compare the model prediction solution, denoted as $u_h$, with a reference solution, denoted as $u_{ref}$. The reference solution can be obtained from an exact solution or numerically computed via methods such as finite difference or spectral methods. For each numerical experiment, we present the following three types of errors for the predicted solutions:
  \beqs
  && \text{Relative~} L^2 \text{~error} := \displaystyle \frac{\|u_h-u_{ref}\|_2}{\|u_{ref}\|_2}; \\
  && \text{Relative~} L^1 \text{~error} := \displaystyle \frac{\|u_h-u_{ref}\|_1}{\|u_{ref}\|_1};\\
  && \text{Max error} := \|u_h-u_{ref}\|_{\infty}.
  \eeqs

  We also conducte robustness tests for several numerical examples. The method we propose not only achieves satisfied prediction results on the target PDE, but also demonstrates good prediction accuracy for functions near the target function, owing to the generalization ability of the neural network. To achieve this, we introduce noise to the target initial value or target source term to create new test functions. Specifically, we select $n$ pairs of data $(x_i, y_i)_{i=1,2,\cdots, n}$, and then add noise to $y_i$ to obtain $y_i^{new} = y_i+\gamma*\mathcal{N}(0,\sigma^2_y)$, where $\sigma_y^2$ is the variance of the entire dataset. Subsequently, we interpolate $(x_i, y_i^{new})$ into a smooth function using cubic interpolation as the new test function. We test the accuracy for these new test functions across several PDE examples, where the corresponding reference solution can be computed using the pseudospectral method or the artificial boundary method.

  In addition to demonstrating the capability of the proposed method to solve various equations on unbounded domains, we also investigate the impact of important factors such as the size of the spatial domain where the data is generated, the number of input functions, and physics-informed training on the experimental results. The findings and results of these investigations are provided in Appendix B, C, and D, respectively.
  
  \subsection{1-D wave equation}
  We consider the 1-D second-order wave equation defined on the whole space:
  \begin{eqnarray}\label{wave}
  \begin{aligned}
  & u_{tt}-u_{xx}=f(x,t), && x\in \mathbb{R},~ t\in(0,1], \\
  & u(x, 0) = \varphi_0(x), && x\in \mathbb{R},\\
  & u_t(x,0) = \varphi_{1}(x), && x\in \mathbb{R}.
  \end{aligned}
  \end{eqnarray}
  We solve two representative examples: one with a known analytical solution, and the other without an analytical solution. For the latter, we compute a reference solution using the finite difference method on a sufficiently large domain. To generate training data for both examples, we construct various forms of analytical expressions.
  
  \subsubsection{Case 1: The analytical solution exists}
\begin{figure}[H]
\begin{center}
  \begin{minipage}[c]{0.55\textwidth}
    \centering
    \includegraphics[width=0.48\textwidth]{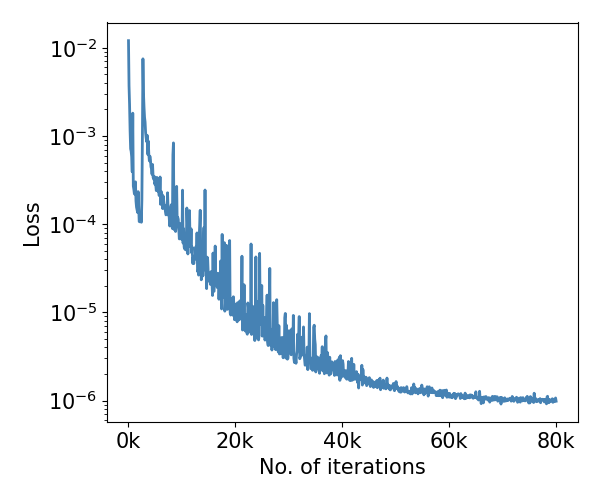}
    \includegraphics[width=0.48\textwidth]{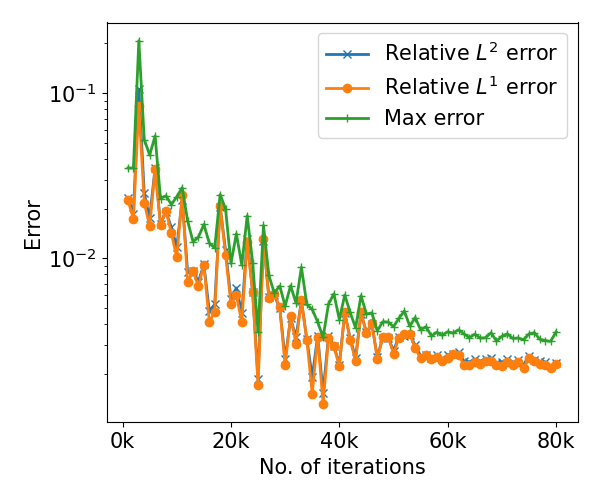}
\end{minipage}%
\hfill
\begin{minipage}[c]{0.44\textwidth}
  \footnotesize
    \centering
    \begin{tabular}{cccc}
        \toprule
         & $L^2$ error & $L^1$ error  & max error\\
        \midrule
        target PDE & 2.34e-03 & 2.30e-03 & 3.60e-03  \\
        $1\%$ noise on $f$ & 2.53e-03 & 2.43e-03  & 4.40e-03 \\
        $10\%$ noise on $f$ & 8.96e-03 & 7.25e-03  & 1.92e-02 \\
        $1\%$ noise on $\varphi_1$ & 2.59e-03 & 2.49e-03  & 4.58e-03 \\
        $10\%$ noise on $\varphi_1$ & 1.04e-02 & 8.49e-03  & 2.24e-02 \\
        \bottomrule
    \end{tabular}
\end{minipage}
\end{center}
\caption{\textit{Solving a 1D wave equation (case 1).} Left: Training loss and test errors for 80,000 iterations. Right: Test errors on the target PDE w/ and w/o noise.}
\label{wave1d-1-loss}
\end{figure}

\begin{figure}[H]
    \centering
\includegraphics[width=1.0\linewidth]{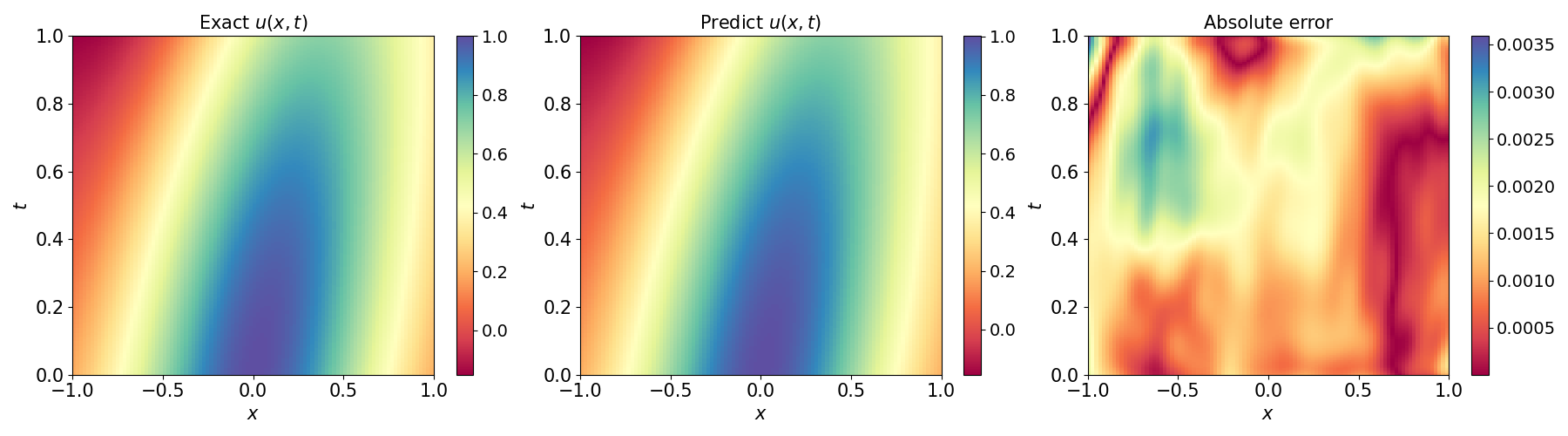}
\caption{\textit{Solving a 1D wave equation (case 1):} Exact solution versus the predictions of the trained MIONet.}
\label{wave1d-1-sol}
\end{figure}

  We consider the following initial conditions and the source term:
  \beq\label{wave_case1_target}
  \begin{aligned}
  &\varphi_0(x) = \exp\left(-x^2\right)\cos(x),\\
  &\varphi_{1}(x) = \exp\left(-x^2\right)\sin(x),\\
  &f(x,t) = \exp(-x^2)\left(4x\sin(t-x)+(2-4x^2)\cos(t-x)\right),
  \end{aligned}
  \eeq
  with the exact solution $u(x,t) = \exp\left(-x^2\right)\cos(t-x)$. The computational domain is $(x,t)\in [-1,1]\times[0,1]$.
  We choose the following form of analytical solutions to generate the required data:
  \begin{eqnarray}\label{wave_case1_construct_sol}
  u(x,t) = \sum_{i=0}^K A_i\sin(k_it+a_i)H_i(wx+b),
  \end{eqnarray}
  where $H_i(x)$ stands for the Hermite funtion \cite{shen2011spectral}, $A_i, k_i, a_i, w, b$ are parameters randomly generated. We set $K=2$. The spatial domain of generating data is $[-3,3]$. We select 5000 input function pairs, i.e., $N=5000$. The training loss and three-type errors computed at a $101\times201$ spatio-temporal grid are presented in Figure \ref{wave1d-1-loss}, with the error computed at every 1,000 iterations. It is observed that as the loss stabilizes during training, the testing error also tends to stabilize. This implies that monitoring the loss can provide a rough indication of whether the trained model is suitable for the target problem. The average values of the last 20 relative $L^2$ errors, relative $L^1$ errors and max errors shown in Figure \ref{wave1d-1-loss} are 2.34e-03, 2.30e-03, and 3.60e-03, respectively. 
  The results of the robustness test are depicted on the right side of Figure \ref{wave1d-1-loss}. We evaluate the prediction accuracy of the model by introducing 1\% and 10\% noise to source term $f$ and initial value $\varphi_1$, respectively. The errors shown are the average of the errors from 100 different noisy test functions.
  From the data, it is evident that when 1\% noise is added, the test error remains relatively close to that of the target PDE. However, with the introduction of 10\% noise, the test error increases, albeit maintaining a certain level of accuracy.
  Figure \ref{wave1d-1-sol} illustrates visual representations of the exact solution, the prediction solution and the absolute error for target PDE.

  \subsubsection{Case 2: Point source and without analytical solution}
\begin{figure}[H]
  \begin{center}
    \begin{minipage}[c]{0.54\textwidth}
      \centering
      \includegraphics[width=0.48\textwidth]{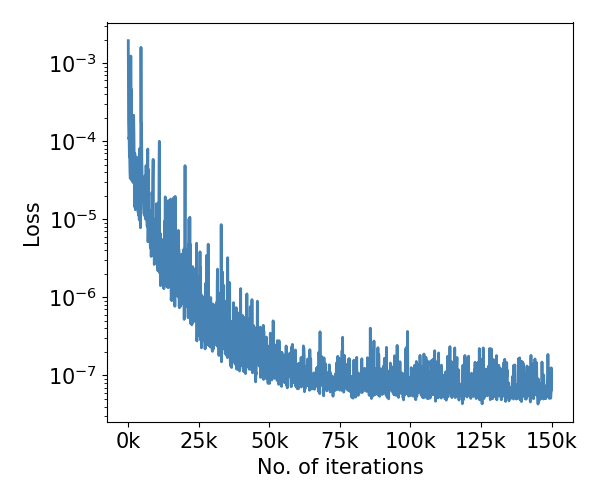}
      \includegraphics[width=0.48\textwidth]{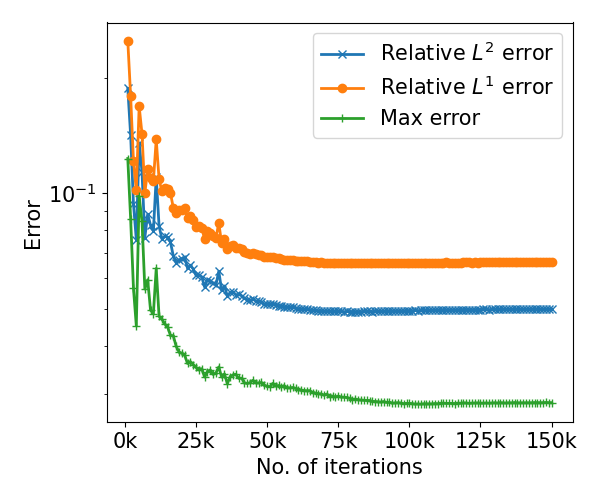}
  \end{minipage}%
  \hfill
  \begin{minipage}[c]{0.44\textwidth}
    \footnotesize
      \centering
      \begin{tabular}{cccc}
          \toprule
           &  $L^2$ error & $L^1$ error  & max error\\
          \midrule
          target PDE & 4.98e-02 &6.60e-02 & 2.84e-02  \\
          $1\%$ noise on $f$ & 4.69e-02 & 6.27e-02  & 2.66e-02 \\
          $10\%$ noise on $f$ & 5.01e-02 & 6.54e-02  & 3.30e-02 \\
          \bottomrule
      \end{tabular}
  \end{minipage}
  \end{center}
  \caption{\textit{Solving a 1D wave equation (case 2).} Left: Training loss and test errors for 150,000 iterations. Right: Test errors on the target PDE w/ and w/o noise.}
  \label{wave1d-2-loss}
\end{figure}

\begin{figure}[H]
      \centering
  \includegraphics[width=1.0\linewidth]{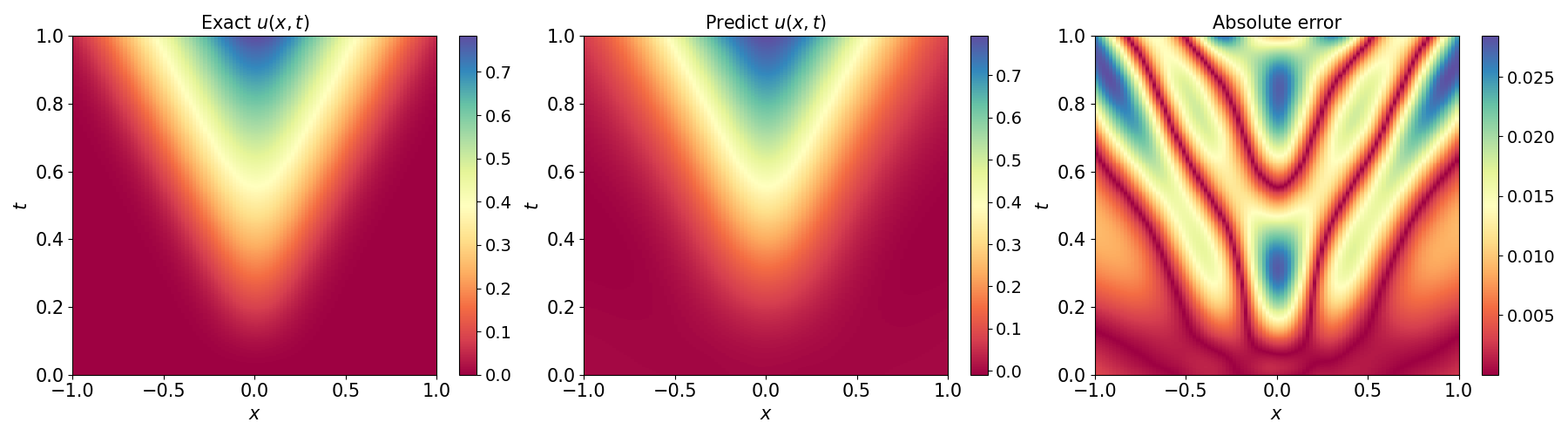}

  \caption{\textit{Solving a 1D wave equation (case 2):} Exact solution versus the predictions of the trained MIONet.}
  \label{wave1d-2-sol}
  \end{figure}

  We consider zero initial values, i.e., $\varphi_0(x) = 0, \varphi_1(x) = 0$, and $f(x,t) = 5\exp(-25x^2)$, which is an approximation of point source function.
  The reference solution $u(x,t)$ of this case is computed using a second-order finite difference method on a large enough spatial domain. To generate the approximate source term $f$, we choose the following form of the analytical solutions: 
  \begin{eqnarray*}
  u(x,t) = \sum_{i=0}^K A_it^2\cos(a_it+b_i)\exp(-x^2/\sigma_{i}^2).
  \end{eqnarray*}
  The spatial domain of generating data is $[-2,2]$, and $K=2$. We select 5000 input function pairs, i.e., $N=5000$. The training loss and the error computed at a $101\times201$ spatio-temporal grid are presented in Figure \ref{wave1d-2-loss}. The average values of the final 20 errors are 4.98e-02, 6.60e-02, and 2.84e-02 for the relative $L^2$ error, relative $L^1$ error, and max error, respectively. 
  Compared to case 1, it is relatively difficult to generate the appropriate analytical expression of this case due to the special form of the function $f$. As a result, the prediction error is also greater. To achieve higher accuracy, it is crucial to construct analytical expressions of higher quality. The results of the robustness test are shown on the right side of Figure \ref{wave1d-2-loss}. It can be observed that when adding 1\% and 10\% noise to f, the prediction error of the model remains at the same level.
  Finally, we draw the reference solution, the prediction solution, and the absolute error of the target PDE in Figure \ref{wave1d-2-sol}.

\subsubsection{Case 3: Solution with a large compact support}

\begin{figure}[ht]
\centering
\includegraphics[width=0.4\linewidth]{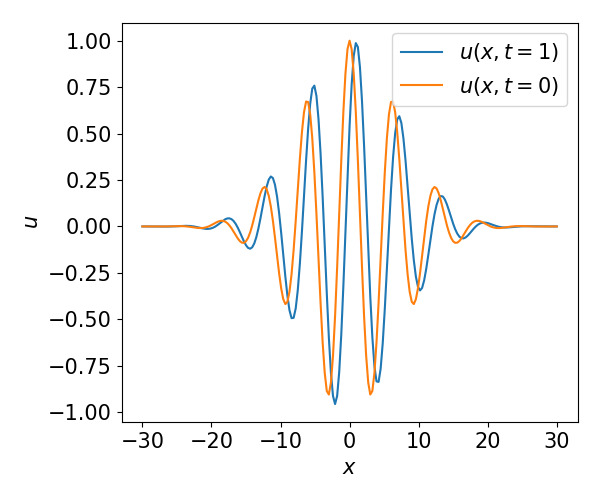}
\caption{\textit{Solving a 1D wave equation (case 3):} Exact solution.}
\label{wave1d-3-u}
\end{figure}

Artificial boundary methods typically require that the support sets of initial values and source terms be contained within the computational domain, posing challenges for problems with large support sets. However, our proposed method in this work overcomes this limitation. We illustrate its applicability to such problems by examining a case with a large compact support set. Concretely, we consider the wave equation \eqref{wave} with the following analytical solution:
\beq
u(x,t) = \exp(-0.01x^2)\cos(t-x),
\eeq
When $|x|\approx40$, the value of the solution only decays to $10^{-6}$. Figure \ref{wave1d-3-u} visually represents the initial data and the state of $u(x,t=1)$. 
  
\begin{figure}[ht]
\centering
\includegraphics[width=0.45\linewidth]{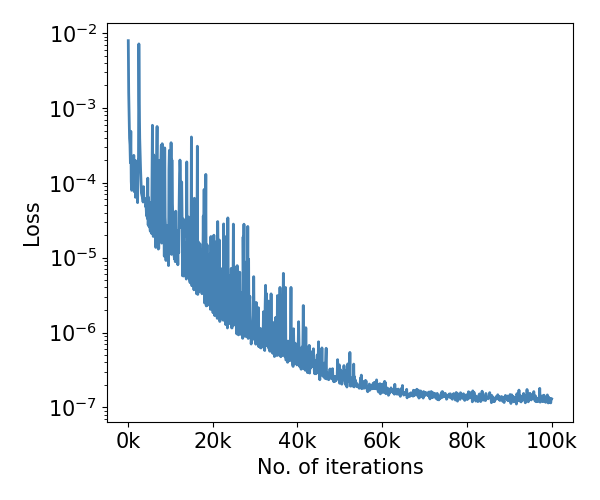}
\includegraphics[width=0.45\linewidth]{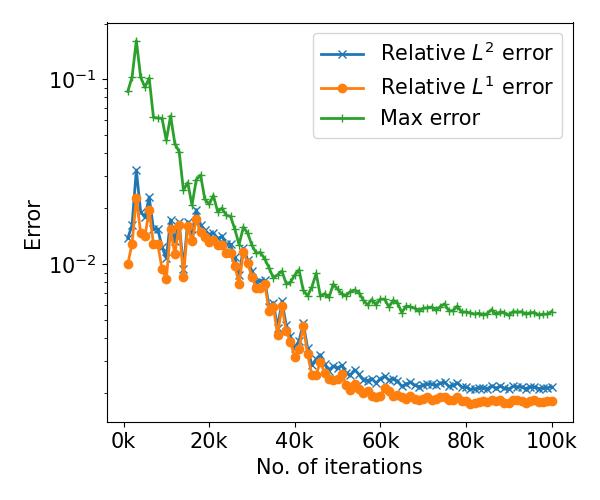}
\caption{\textit{Solving a 1D wave equation (case 3):} Training loss and test errors for 100,000 iterations.}
\label{wave1d-3-loss}
\end{figure}

\begin{figure}[ht]
\centering
\includegraphics[width=1.0\linewidth]{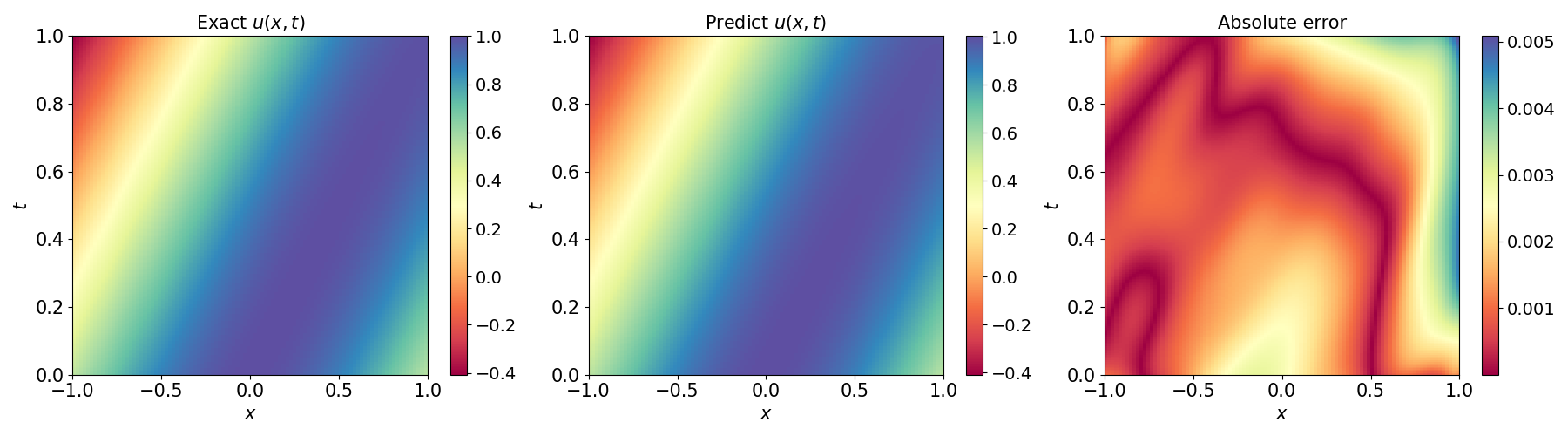}
\caption{\textit{Solving a 1D wave equation (case 3):} Exact solution versus the predictions of the trained MIONet.}
\label{wave1d-3-sol}
\end{figure}

The domain of interest we consider is $[-1,1]\times[0,1]$. We generate 1000 training data from the following expression
\beqs
u(x,t) = \sum_{i=0}^{K}A_i\exp(a_i(x+c_it)^2)\sin(k_ix+b_it+d_i).
\eeqs
The training loss and errors computed at a $101\times201$ spatio-temporal grid are presented in Figure \ref{wave1d-3-loss}. Figure \ref{wave1d-3-sol} displays the reference solution, prediction solution and the absolute error at final time $t=1$. Furthermore, we calculated the average values of the final 20 relative $L^2$ errors, relative $L^1$ errors, and max errors as 2.17e-03, 1.82e-03, and 5.54e-03, respectively. 


\subsubsection{Case 4: Non decaying initial values}
\begin{figure}[H]
  \begin{center}
    \begin{minipage}[c]{0.55\textwidth}
      \centering
      \includegraphics[width=0.48\textwidth]{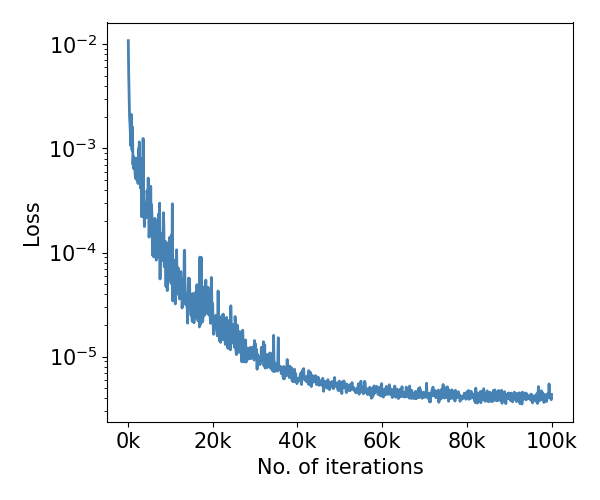}
      \includegraphics[width=0.48\textwidth]{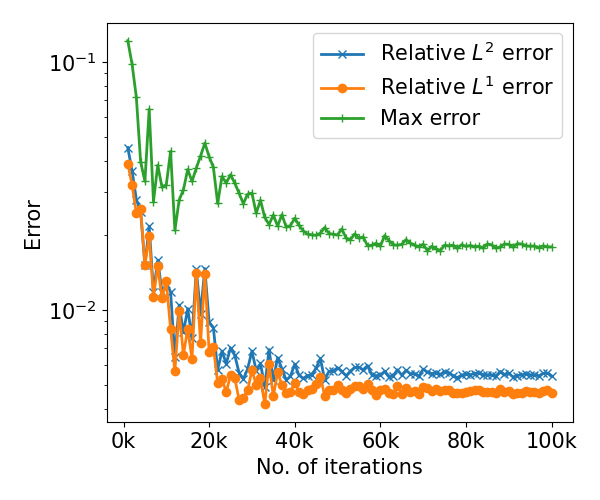}
  \end{minipage}%
  \hfill
  \begin{minipage}[c]{0.44\textwidth}
    \footnotesize
      \centering
      \begin{tabular}{cccc}
          \toprule
           & $L^2$ error & $L^1$ error  & max error\\
          \midrule
          target PDE & 4.73e-03 & 4.16e-03 & 2.07e-02  \\
          $1\%$ noise on $f$ & 5.59e-03 & 4.74e-03  & 1.82e-02 \\
          $10\%$ noise on $f$ & 6.31e-03 & 5.32e-03  & 1.93e-02 \\
          $1\%$ noise on $\varphi_1$ & 6.33e-03 & 5.00e-03  & 2.33e-02 \\
          $10\%$ noise on $\varphi_1$ & 7.20e-02 & 4.76e-02  & 3.26e-01 \\
          \bottomrule
      \end{tabular}
  \end{minipage}
  \end{center}
  \caption{\textit{Solving a 1D wave equation (case 4).} Left: Training loss and test errors for 100,000 iterations. Right: Test errors on the target PDE w/ and w/o noise.}
  \label{wave1d-4-loss}
\end{figure}

\begin{figure}[ht]
  \centering
  \includegraphics[width=1.0\linewidth]{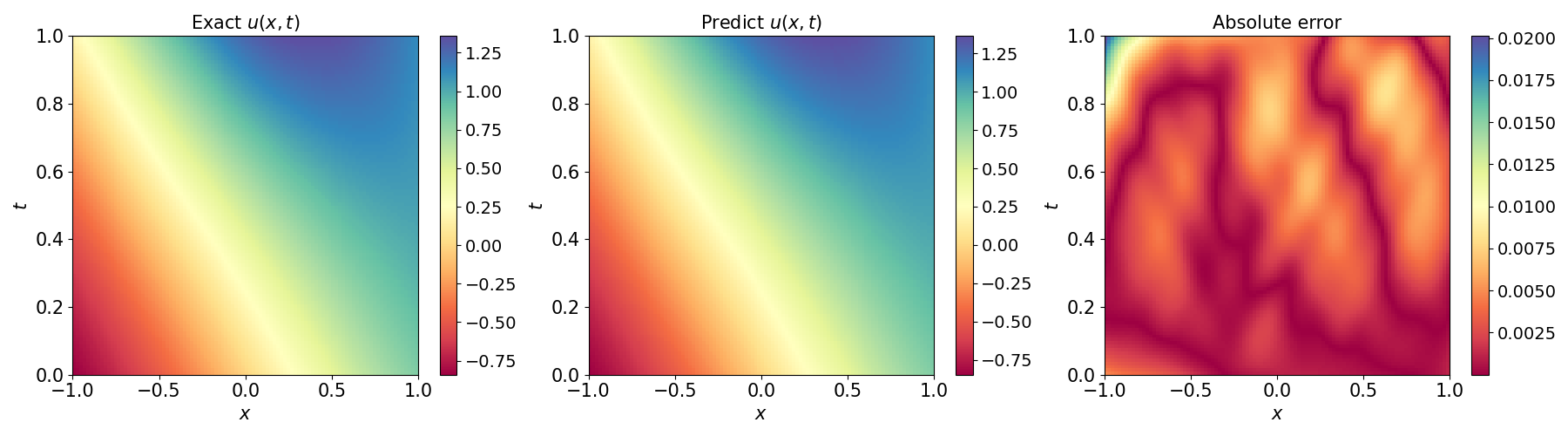}
  \caption{\textit{Solving a 1D wave equation (case 4):} Exact solution versus the predictions of the trained MIONet.}
  \label{wave1d-4-sol}
  \end{figure}

We consider a case where initial values do not decay at infinity:
\beq\label{wave_case4_target}
\begin{aligned}
&\varphi_0(x) = \sin(x),\\
&\varphi_{1}(x) = \cos(x),\\
&f(x,t) = \exp(-x^2).
\end{aligned}
\eeq
We utilize the artificial boundary method \cite{han2013artificial} to compute the reference solution.
The domain of interest is $[-1,1]\times[0,1]$. We generate 2000 training data from the following expression
\beqs
u(x,t) = A \sin(k_1 x+w_1t+b_1)+\cos(k_2x+w_2t+b_2).
\eeqs
The training loss and errors computed at a $101\times201$ spatio-temporal grid are presented in the left of Figure \ref{wave1d-4-loss}. The average values of the final 20 relative $L^2$ errors, relative $L^1$ errors, and max errors as 4.73e-03, 4.16e-03, and 2.07e-02, respectively. The robustness test results are shown in the right of Figure \ref{wave1d-4-loss}.
Figure \ref{wave1d-4-sol} displays the reference solution, prediction solution and the absolute error of the target PDE. 

\subsection{2D wave equation}
  This example aims to solve the following 2D wave equation 
      \begin{eqnarray*}\label{wave2d}
      \begin{aligned}
      & u_{tt}-\Delta u=f(\vx,t), && \vx\in \mathbb{R}^2,~ t\in(0,1], \\
      & u(\vx, 0) = \varphi_0(\vx), && \vx\in \mathbb{R}^2,\\
      & u_t(\vx,0) = \varphi_{1}(\vx), && \vx\in \mathbb{R}^2.
      \end{aligned}
      \end{eqnarray*}
  The target initial functions and source are given as:
  \beqs
  && \varphi_0(\vx) = \exp\left(-\frac{x_1^2+x_2^2}{2}\right),\\
  && \varphi_1(\vx) = \exp\left(-\frac{x_1^2+x_2^2}{2}\right)(x_1+x_2),\\
  && f(\vx, t) = \exp\left(-\frac{(x_1-t)^2+(x_2-t)^2}{2}\right)\left((x_1-t+x_2-t)^2\cos(kt)-2k(x_1-t+x_2-t)\sin(kt)\right. \\
  && \quad\quad\quad\quad \left.-(k^2+(x_1-t)^2+(x_2-t)^2)\cos(kt)\right),
  \eeqs
  which are computed from the solution $$u(\vx,t)=\exp\left(-\frac{(x_1-t)^2+(x_2-t)^2}{2}\right)\cos(kt).$$
  The spatial domain of interest is selected as $[-1,1]$. 
  Based on the property of the target function, we construct the following form of analytical solutions to generate training data:
  \beqs
  u(\vx,t) = A\exp\left(-\frac{(x_1-a_1t)^2+(x_2-a_2t)^2}{\sigma^2}\right)\cos(k_1x_1+k_2x_2-wt).
  \eeqs
  
\begin{figure}[H]
      \centering
      \includegraphics[width=0.45\linewidth]{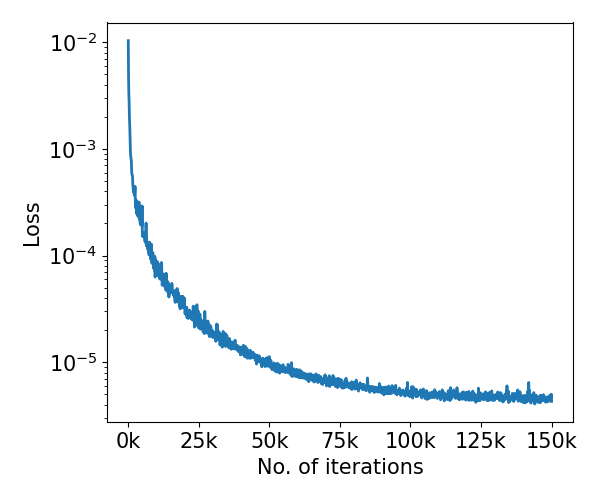}
      \includegraphics[width=0.45\linewidth]{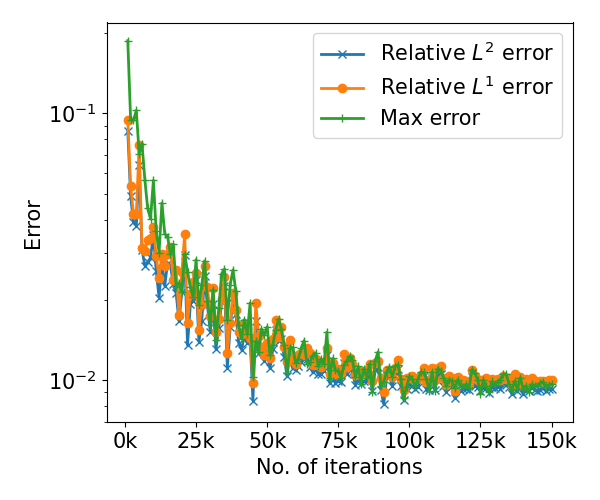}
      \caption{\textit{Solving a 2D wave equation:} Training loss and test errors for 150,000 iterations.}
      \label{wave2d-loss}
      \end{figure}

\begin{figure}[H]
          \centering
      \includegraphics[width=1.0\linewidth]{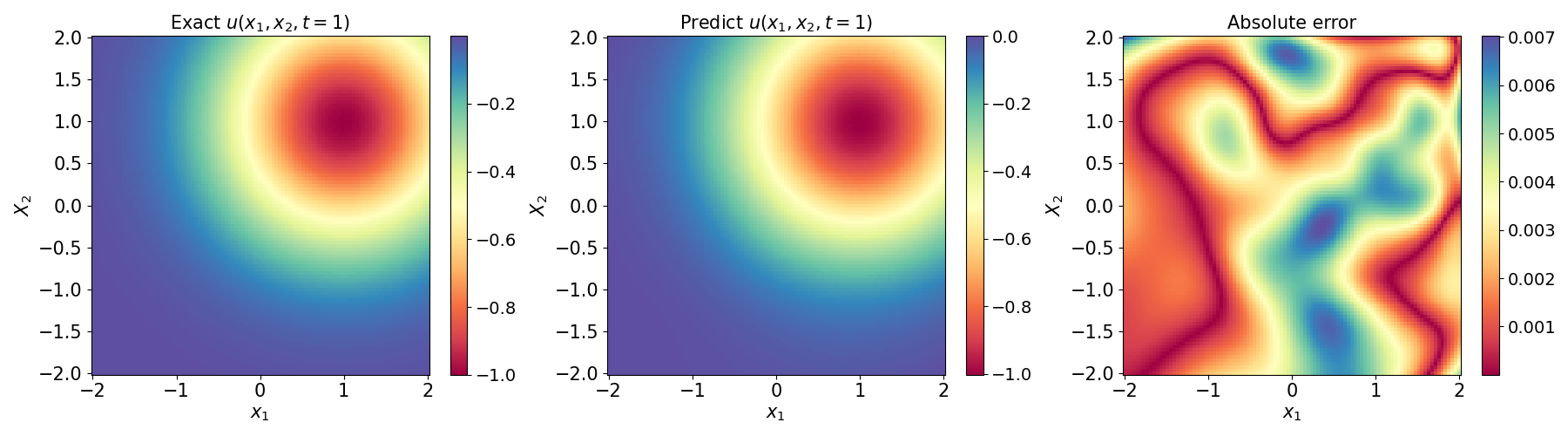}
      \caption{\textit{Solving a 2D wave equation:} Exact solution versus the predictions of the trained MIONet at final time $t=1$.}
      \label{wave2d-sol}
  \end{figure}

  We generate 2000 input function pairs on the spatial domain $[-5,5]\times[-5,5]$. The training loss and errors computed at a $21\times21\times51$ spatio-temporal grid are presented in Figure \ref{wave2d-loss}. Figure \ref{wave2d-sol} displays the reference solution, prediction solution and the absolute error at final time $t=1$. Furthermore, we calculated the average values of the final 20 relative $L^2$ errors, relative $L^1$ errors, and max errors as 9.27e-03, 1.00e-02, and 9.65e-03, respectively.

  \subsection{Burgers' equation}
  To highlight the ability of the proposed method to handle nonlinear equations on unbounded domains, we solve the following 1D Burgers' equation
  \begin{eqnarray}\label{Burgers}
  \begin{aligned}
  & u_{t}+uu_{x}-\nu u_{xx} = f(x,t), && x\in \mathbb{R},~ t\in(0,1], \\
  & u(x, 0) = \varphi(x), && x\in \mathbb{R}.
  \end{aligned}
  \end{eqnarray}
  In this section, we present two examples: solving a single equation and solving a family of parameterized equations. The latter case means that we solve multiple PDEs simultaneously using a single trained model.
  
  \subsubsection{Case 1: Solve one single equation} \label{sec_Burgers1}
\begin{figure}[H]
  \begin{center}
    \begin{minipage}[c]{0.54\textwidth}
      \centering
      \includegraphics[width=0.49\textwidth]{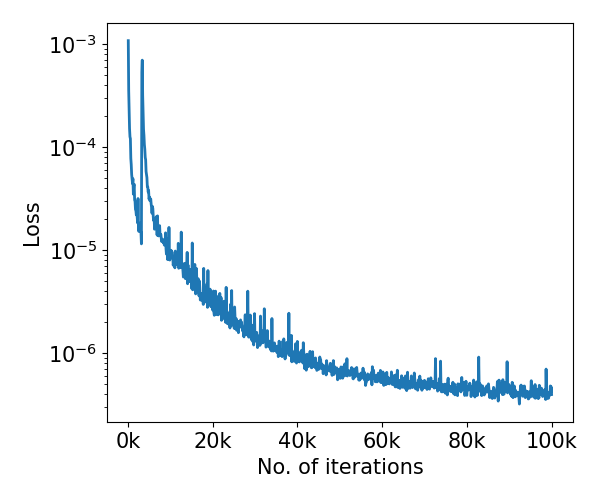}
      \includegraphics[width=0.49\textwidth]{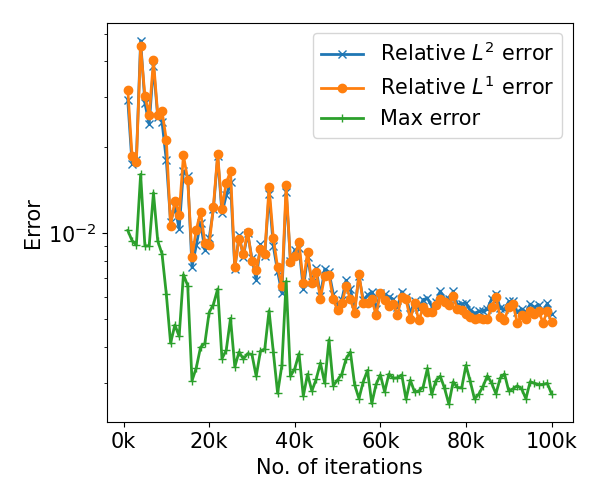}
  \end{minipage}%
  \hfill
  \begin{minipage}[c]{0.44\textwidth}
    \footnotesize
      \centering
      \begin{tabular}{cccc}
          \toprule
           &  $L^2$ error & $L^1$ error  & max error\\
          \midrule
          target PDE & 5.51e-03 & 5.20e-03  & 2.91e-03  \\
          $1\%$ noise on $f$& 6.99e-03 & 6.71e-03  & 3.96e-03 \\
          $10\%$ noise on $f$& 4.79e-02 & 4.89e-02  & 1.82e-02 \\
          \bottomrule
      \end{tabular}
  \end{minipage}
  \end{center}
  \caption{\textit{Solving a Burgers' equation.} Left: Training loss and test errors for 100,000 iterations. Right: Test errors on the target PDE w/ and w/o noise.}
  \label{Burgers1-loss}
\end{figure}

\begin{figure}[H]
\includegraphics[width=1.0\linewidth]{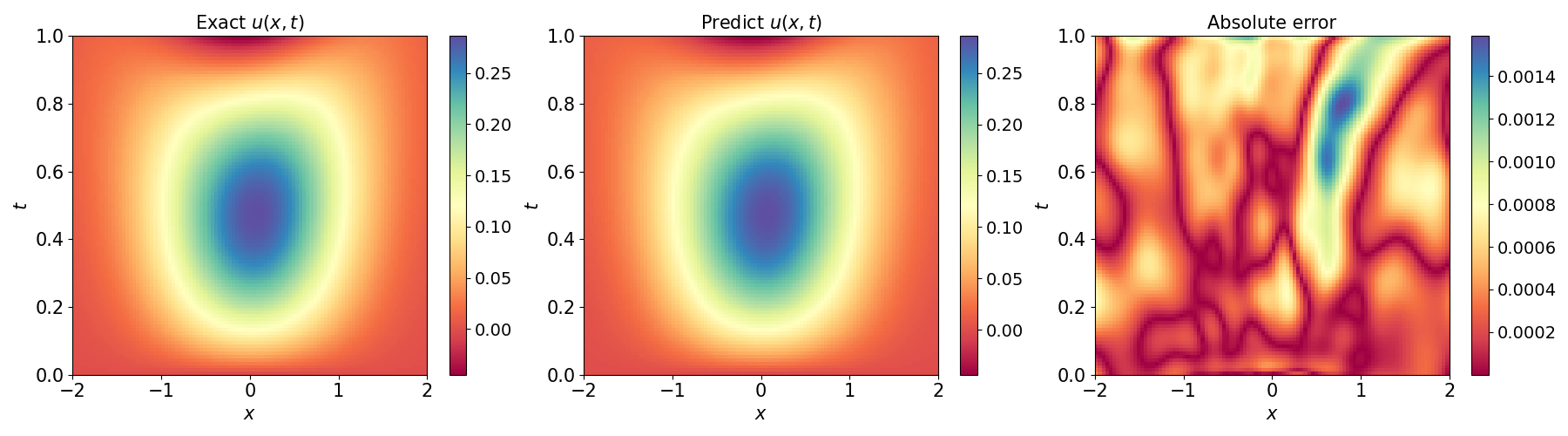}
\caption{\textit{Solving a Burgers' equation:} Exact solution versus the predictions of the trained MIONet.}
\label{Burgers1-sol}
\end{figure}

  In this example, we consider $\varphi(x) = 0,~ f(x,t) = \cos(\pi t)\exp(-x^2)$, the viscosity is set to $\nu = 0.2$. To evaluate the accuracy of the predicted solution, we solve the equation using the finite difference method on a large enough domain, where the spatial step size and temporal step size are both set to 0.01. The domain of interest is selected as $[-2,2]$. We construct the analytical solution to generate training data:
  \beq \label{BurgersData}
  u(x,t) = \sum_{i=0}^K A_i(\sin(k_{1,i}\pi t)+t\cos(k_{2,i}\pi t))H_{i}(w_ix+b_i),
  \eeq
  which spontaneously satisfies the target initial value. Therefore, our focus is on adjusting the parameters of this constructed solution to make it closely match the target source term.
  The spatial domain for generating data is also set as $[-2,2]$. We select 5000 input function pairs, i.e., $N=5000$. We discuss the impact of the number of input functions on model training and prediction accuracy in Appendix C. Generally, increasing $N$ tends to improve the prediction accuracy of the model.
  The training loss and test errors for the target PDE and noisy PDEs are presented in Figure \ref{Burgers1-loss}. Additionally, Figure \ref{Burgers1-sol} illustrates the comparison between the predicted and the exact solution for the target PDE. 
  These results demonstrate that the proposed method can effectively tackle nonlinear equations on unbounded domains, which are often challenging for traditional numerical methods to solve.

  \subsubsection{Case 2: Solve multiple equations simultaneously}
\begin{figure}[ht]
      \centering
      \includegraphics[width=0.32\linewidth]{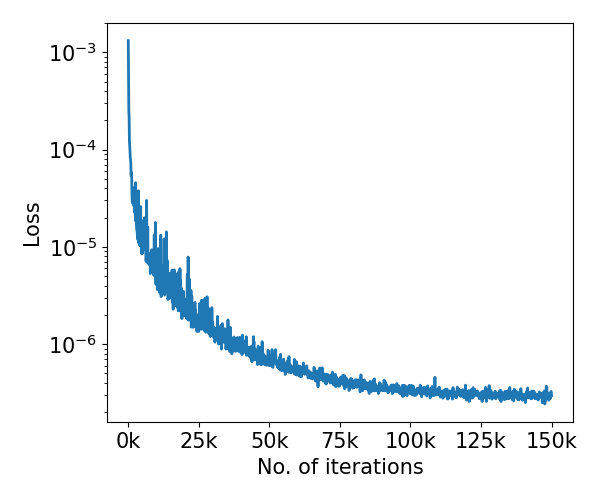}
      \includegraphics[width=0.32\linewidth]{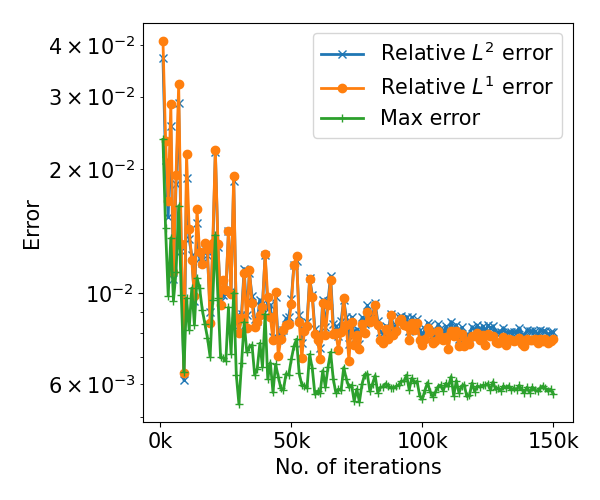}
      \includegraphics[width=0.32\linewidth]{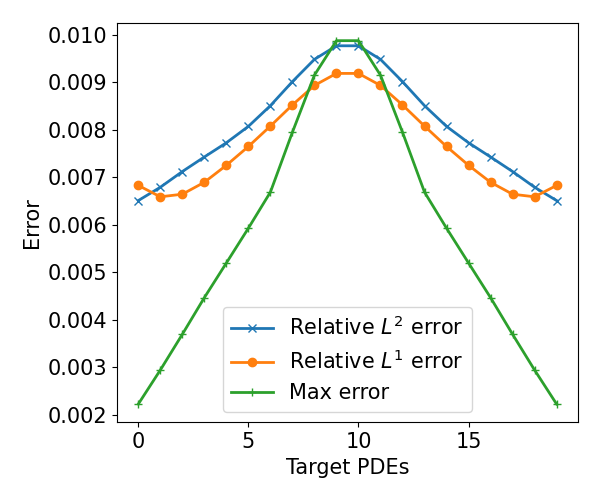}
      \caption{\textit{Solving multiple Burgers' equations:} Left: Training loss for 150,000 iterations. Middle: Average error of predicted solutions for all 20 target PDEs. Right: Final error of the prediction solution for each target PDE.}
      \label{Burgers2}
  \end{figure}

  In this example, we show the proposed method is also capable of solving multiple PDEs simultaneously. We consider Burgers' equation \eqref{Burgers} with a zero initial value and a family of source terms parameterized by $\theta$:
  $$f(x,t) = \cos(\theta t)\exp(-x^2).$$
  We uniformly select 20 values of $\theta$ from the interval $[-\pi, \pi]$. The goal is to train only one model to solve these parameterized PDEs.
  The training data is also generated from the analytical solutions given in \eqref{BurgersData}. The selected solutions satisfy the condition that the error between the corresponding source term and each target source term is below a given tolerance $\epsilon$. The data is generated over the spatial domain $[-4,4]$. In the left plane of Figure \ref{Burgers2}, we plot the training loss over 150,000 iterations.  The middle panel displays the average error of the predicted solutions for all 20 target PDEs. After the model convergence, the average values of the relative $L^2$ errors, relative $L^1$ errors, and max errors are measured as 8.04e-03, 7.64e-03, and 5.83e-03, respectively.
  Lastly, the right plane of Figure \ref{Burgers2} shows the final error of the prediction solution for each target PDE. This example demonstrates that the proposed method still performs well in solving multiple similar PDEs simultaneously.

  \subsection{KdV equation}\label{sec_KdV}
  This example considers the KdV equation, which describes the propagation of waves in certain nonlinear dispersive media.  It takes the following form:
  \begin{eqnarray}\label{KdV}
  \begin{aligned}
  & u_{t}+ 6uu_{x}+u_{xxx} = f(x,t), && x\in \mathbb{R},~ t\in(0,1], \\
  & u(x, 0) = \varphi(x), && x\in \mathbb{R}.
  \end{aligned}
  \end{eqnarray}
  Here the target initial value and source term are generated from the exact solution $u(x,t) = \exp(-(x-t)^2)$.
  Substituting it to the KdV equation \eqref{KdV}, one obtains
  \beq\label{KdV_target}
  \begin{aligned}
  & \varphi(x) = \exp(-x^2),\\
  & f(x,t) = \exp(-(x-t)^2)(12(t-x)\exp(-(t - x)^2)+14(x-t) + 24tx(x-t) + 8(t^3-x^3)).
  \end{aligned}
  \eeq
  
\begin{figure}[H]
    \centering
    \includegraphics[width=0.45\linewidth]{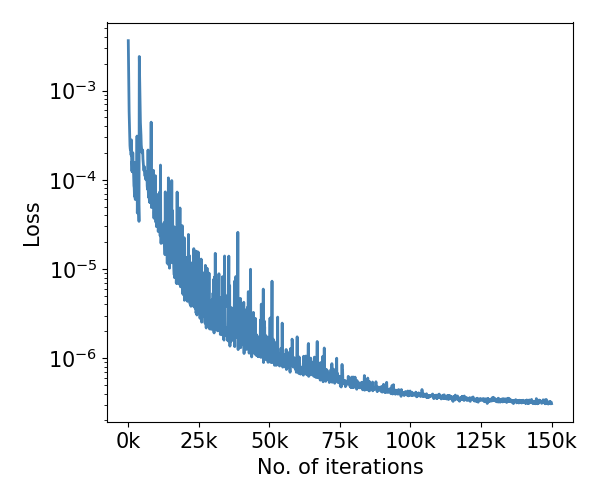}
    \includegraphics[width=0.45\linewidth]{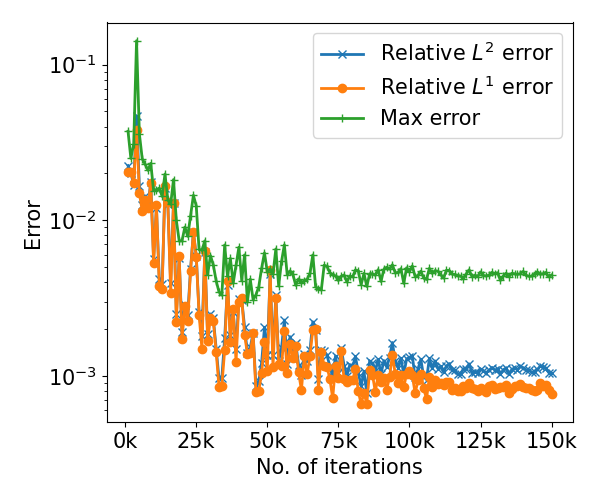}
    \caption{\textit{Solving a KdV equation:} Training loss and test errors for 150,000 iterations.}
    \label{KdV_loss}
  \end{figure}

\begin{figure}[H]
  \includegraphics[width=1.0\linewidth]{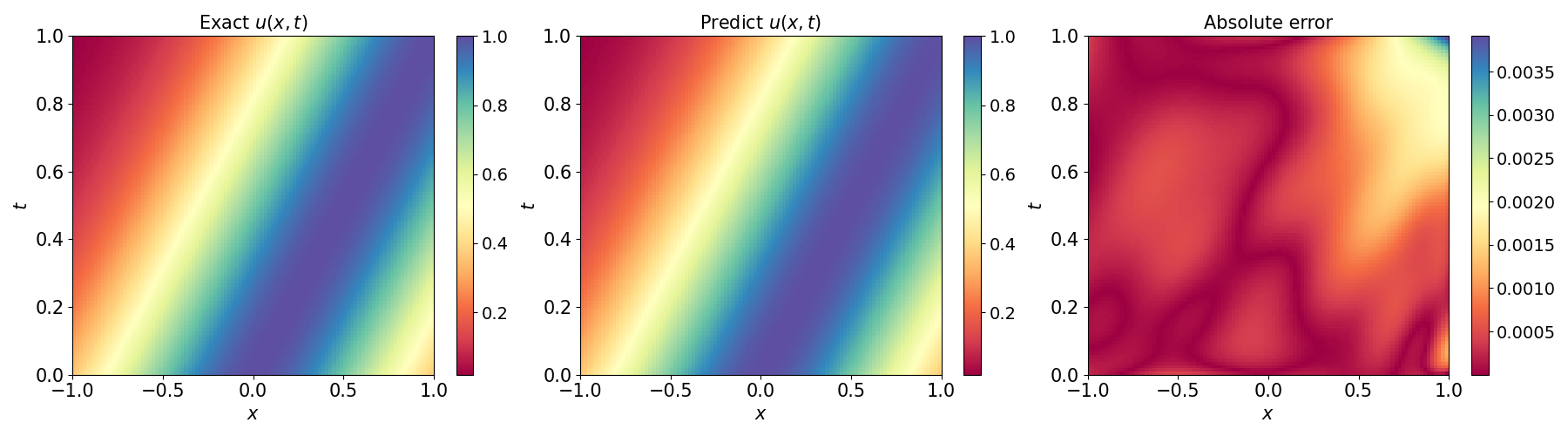}
  \caption{\textit{Solving a KdV equation:} Exact solution versus the predictions of the trained MIONet.}
  \label{KdV_sol}
  \end{figure}

  To solve the equation \eqref{KdV} with the target initial value and source \eqref{KdV_target} on domain $[-1,1]\times[0,1]$, we construct the following analytical expression:
  $$u(x,t)=A\exp(a(x+c_1t+c_2)^2)\cos(kx-wt),$$
  where $A, a, c_1, c_2, k, w$ are tunable parameters. We measure the generated data using the relative $L^2$ error \eqref{metric1} on the spatial domain $[-5,5]$.
  We select $N=1000$ input functions. The training loss and errors tested at a $101\times101$ spatio-temporal grid are presented in Figure \ref{KdV_loss}.  The exact solution, the prediction solution, and the absolute error are drawn in Figure \ref{KdV_sol}. The average values of the relative $L^2$ error, relative $L^1$ error and max error for the last 20 predictions are 1.04e-03, 7.69e-04, and 4.45e-03, respectively.

\begin{figure}[ht]
  \centering
  \includegraphics[width=0.45\linewidth]{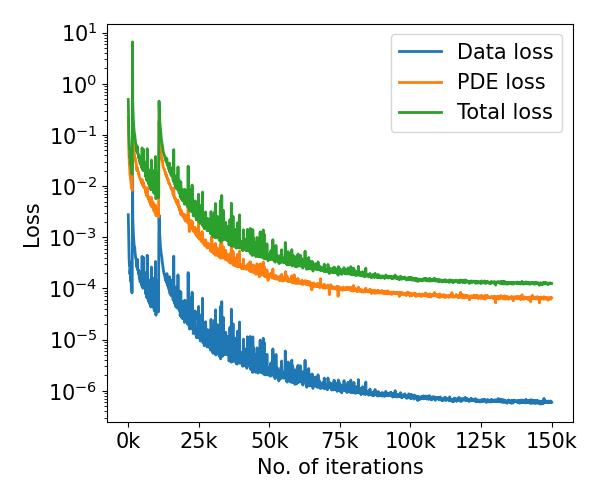}
  \includegraphics[width=0.45\linewidth]{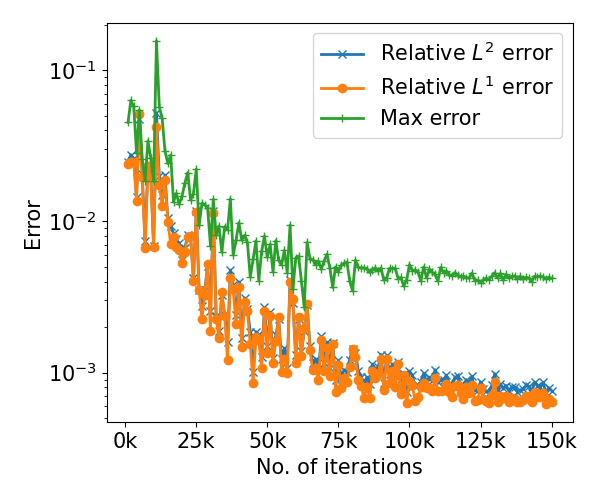}
  \caption{\textit{Solving a KdV equation using the physics-informed training:} Training loss and test errors for 150,000 iterations.}
  \label{KdV_phys_loss}
\end{figure}

\textbf{Physics-informed training.}  
We explore the influence of physics-informed training on prediction results by randomly sampling an additional 500 points to compute the physics loss \eqref{loss_phys}. The model is trained using the total loss \eqref{loss_total}, with weights $w_b$ and $w_d$ set to 1 and 100, respectively. The obtained results are shown in the Figure \ref{KdV_phys_loss}. The average values of the relative $L^2$ error, relative $L^1$ error and max error for the last 20 predictions are 1.03e-03, 8.79e-03, and 5.18e-03, respectively. The error of the predicted solutions obtained by the two training methods is not significantly different.

In practical scenarios where the reference solution is unknown, evaluating accuracy can be challenging. In such cases, we calculate the residual of the solution's equation as a reference. The maximum error of the equation residual for the solution obtained through data training is 4.26, whereas for the solution obtained through physics-informed training, it is 2.76e-02. 
This highlights that training only with data may result in a large equation residual of the predicted solution, even if the error between the prediction and the exact solution is small. This discrepancy arises due to the amplification of errors during the differentiation of the predicted solution. In contrast, results obtained through physics-informed training effectively satisfy the equations, serving as a crucial posterior indicator for predicted solutions. However, it should be noted that the physics-informed training is more time-consuming.

\subsection{Schr\"odinger equation}
\begin{figure}[ht]
      \centering
      \includegraphics[width=1.0\linewidth]{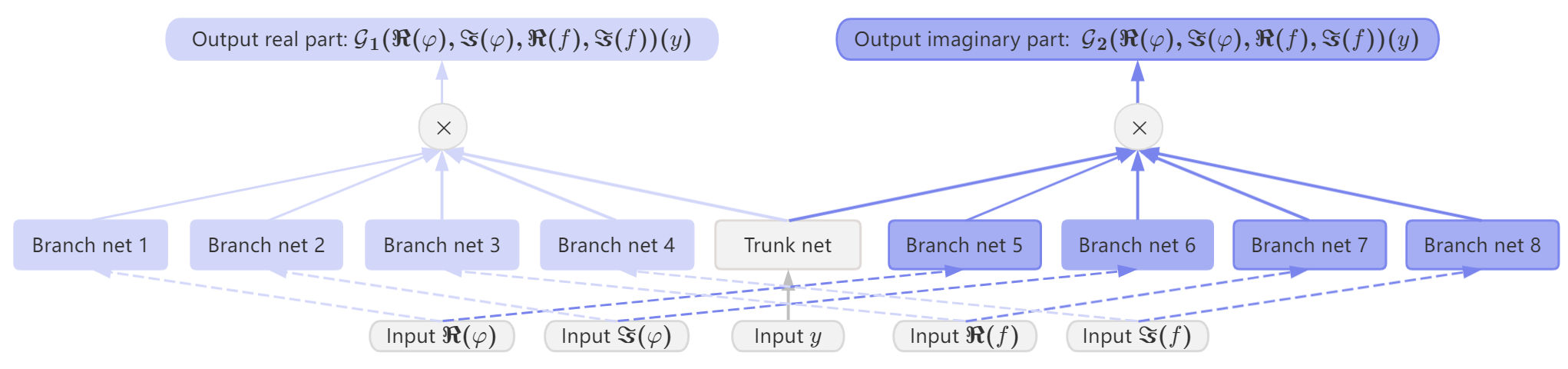}
      \caption{MIONet architecture for solving the Schr\"odinger equation.}
      \label{NN_Schr}
      \end{figure}

  In this example, we consider the  classical Schr\"odinger equation, which is typically defined in an unbounded domain
  \begin{eqnarray}\label{Schr}
  \begin{aligned}
  & {\rm i}u_{t}= -u_{xx} + f(x,t), && x\in \mathbb{R},~ t\in(0,1], \\
  & u(x, 0) = \varphi(x), && x\in \mathbb{R},
  \end{aligned}
  \end{eqnarray}
  where the solution $u$ is the complex-value function. 
  To handle complex values in neural networks, certain modifications need to be made to the network structure since direct gradient propagation is not supported for complex values. The approach involves decomposing complex numbers into their real and imaginary parts, feeding them as inputs to the neural networks, and producing separate outputs for the real and imaginary parts.
  For each input function, there are two corresponding subnetworks: one responsible for outputting the real part and another for outputting the imaginary part. The trunk network encoded input $y$ is shared between these two subnetworks. Figure \ref{NN_Schr} shows in detail the structure of this network.
  
\begin{figure}[ht]
  \centering
  \includegraphics[width=0.45\linewidth]{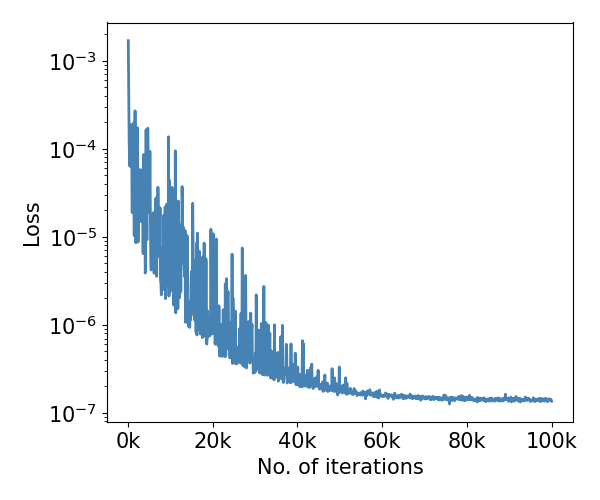}
  \includegraphics[width=0.45\linewidth]{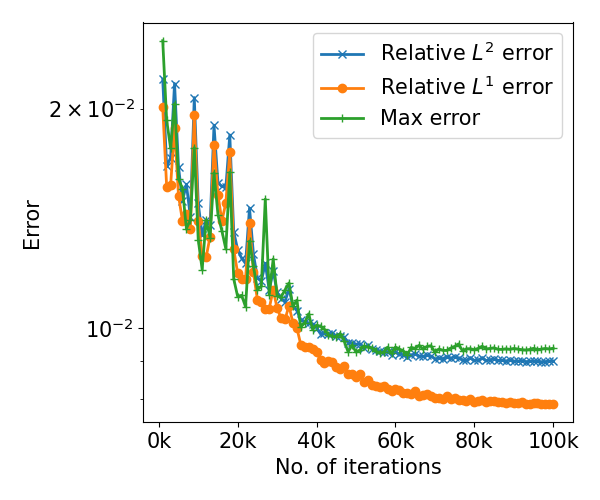}
  \caption{\textit{Solving a 1D Schr\"odinger equation:} Training loss and test errors for 100,000 iterations.}
  \label{Schr-loss}
  \end{figure}

\begin{figure}[ht]
  \centering
  \includegraphics[width=1.0\linewidth]{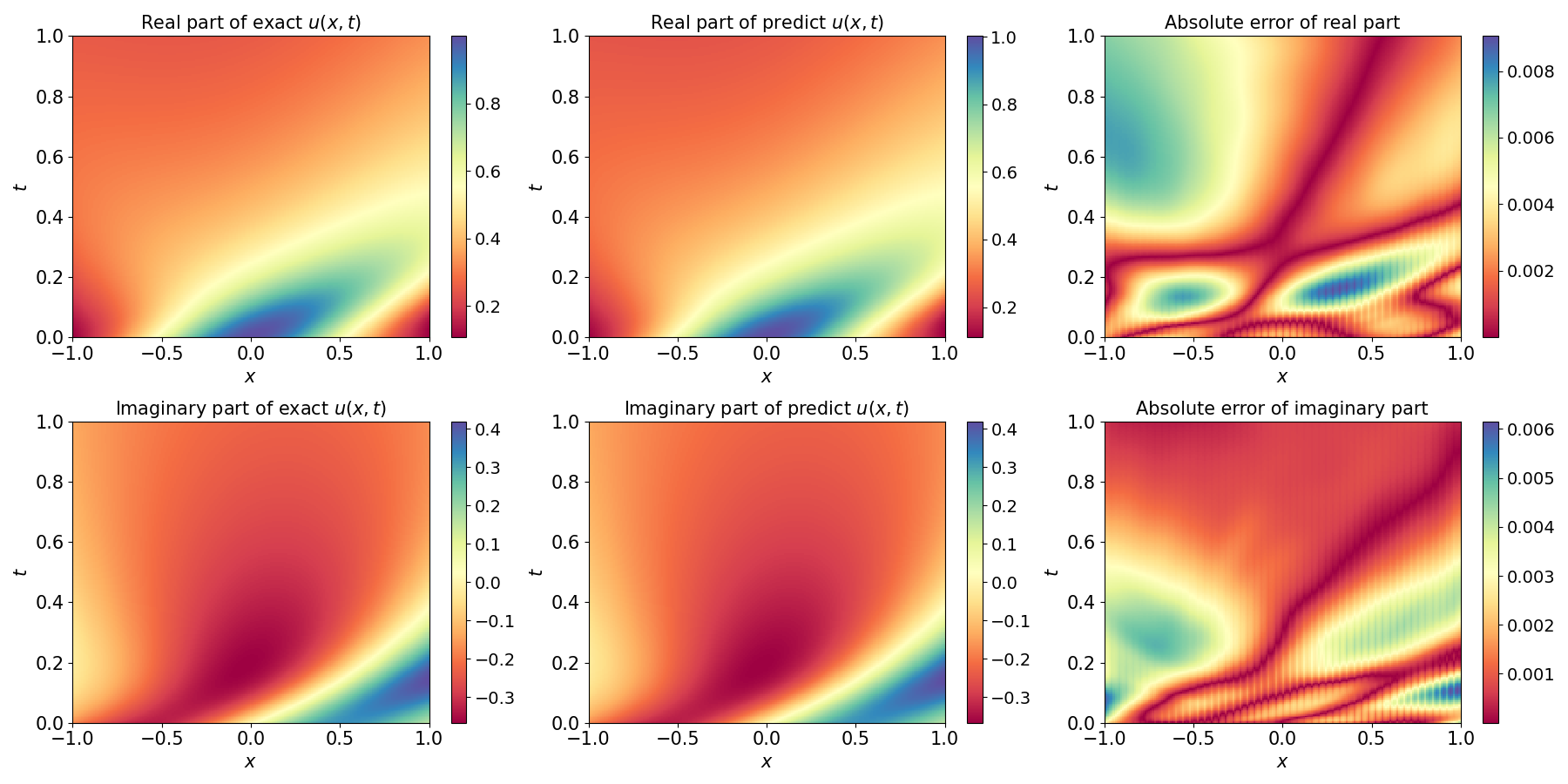}
  \caption{\textit{Sovling a 1D Schr\"odinger equation:} Exact solution versus the predictions of the trained MIONet.}
  \label{Schr-sol}
  \end{figure}

  The target initial value and source term are given as $\varphi(x) = \exp(-x^2+{\rm i}x)\cos(x)$ and $f=0$. 
  It is known that Eq.\eqref{Schr} with $f=0$ has the following analytical solutions:
  \beqs
  \begin{aligned}
  & u(x, t)=\frac{1}{\sqrt{\zeta+{\rm i} t}} \exp \left[{\rm i} k(x-k t)-\frac{(x-2 k t)^2}{4(\zeta+{\rm i} t)}\right],
  \end{aligned}
  \eeqs
  where $k$ is a real parameter that controls the beam propagation speed, and $\zeta$ is a positive parameter that controls the beam width.
  The training data are generated based on this family of solutions, and we remark that the target function is not contained in the training data.
  The spatial domain of generating data is $[-2,2]$. We select 5000 input function pairs, i.e., $N=5000$. The training loss and errors computed at a $101\times101$ spatio-temporal grid in the training process are presented in Figure \ref{Schr-loss}. We evaluate the average values of the final 20 relative $L^2$ errors, relative $L^1$ errors, and max errors, which are found to be 8.90e-03, 7.90e-03, and 9.05e-03, respectively.   
  In the upper part of Figure \ref{Schr-sol}, the real parts of the reference solution and the predicted solution are displayed, accompanied by their respective absolute error. The lower part of Figure \ref{Schr-sol} illustrates the imaginary parts of the reference solution and the predicted solution, along with their corresponding absolute error.

  \section{Conclusion and discussion}
  This paper presents a data generation-based operator learning method for effectively solving PDEs on unbounded domains. 
  The key ingredient of this approach is to generate high-quality training data based on the target PDE. Then, the MIONet model is trained to learn the mapping from the initial value and source term to the solution of the PDE. The generalization ability of the model allows it to directly operate on the initial value and source term of the target PDE, producing a solution with a certain level of accuracy.
  Through extensive testing on different equations, including both linear and nonlinear cases, we have demonstrated the effectiveness of our proposed method.
  
  This study is an attempt to solve PDEs on unbounded domains entirely based on operator learning. Unlike classical numerical methods, this method is not sensitive to the specific form of the equation, making it suitable for solving challenging problems, such as nonlinear cases.
  In addition, the proposed method can not only solve a single PDE but also offers the capability to simultaneously solve multiple PDEs with different parameters. This versatility highlights the wide-ranging applicability and potential impact of our approach.
  Moreover, the cost of data generation in this method is significantly lower compared to the original DeepONet approach which typically relies on expensive traditional numerical methods to generate a large volume of data.  The concept of generating approximate data not only applies to solving problems on unbounded domains but also offers valuable insights for addressing problems on bounded domains and even inverse problems.

Despite the promising results, there are still numerous open questions that require further investigation. 
A key concern is that the proposed method cannot absolutely guarantee the reliability of the predicted solution. Although the physics-informed training enables the evaluation of solution reliability through computing equation residuals, it is crucial to recognize that a solution satisfying the equation may not precisely reflect the original solution on an unbounded domain. 
A key concern is that the proposed method cannot absolutely guarantee the reliability of the predicted solution. The prediction results obtained through the physics-informed training can be substituted into the equation to calculate equation residuals, providing a partial assessment of solution reliability based on the residuals' magnitude. However, this is not a completely reliable evaluation method because a solution  that merely satisfies the equation on an bounded domain may not align with the original equation's solution as defined in an unbounded domain.
Another important concern is that constructing satisfactory analytical solutions is not always straightforward. For instance, in cases where the source term of the target PDE is zero, it is challenging to construct analytical solutions such that the corresponding source terms approximate the target source term. 
In such situations, we must utilize all available information to carefully construct analytical solutions. 
The development of effective construction methods requires further exploration in the future.

\section*{Acknowledgments}
This work is supported by the China Postdoctoral Science Foundation (No. 2023M743266), Zhejiang Provincial Postdoctoral Research Project Merit-based Funding (No.ZJ2023067), Exploratory Research Project (No.2022RC0AN02) and Research Initiation Project (No. K2022RC0PI01) of Zhejiang Lab. 

\bibliographystyle{plain}
\bibliography{DNN,ABM}

  
  \section*{A~~Hyper-parameter settings}
  This appendix provides the detailed parameters used in all examples, including the size of training data, MIONet architectures, learning rate, and the rule of generating training data.
  
  \subsection*{A.1~~Training data.} 
  The number of input function pairs for training is denoted as $N$. We represent the input functions using the function values at $m$ uniform grid points, which we refer to as ``sensors", following the convention established in \cite{lu2021learning}. The training locations for computing initial loss are randomly selected from the spatial domain and locations for computing interior loss are uniformly selected from the spatial-temporal domain. Table \ref{appendix_data} presents the number of input functions $N$, the number of sensors $m$ of input functions (initial function and source function), and the number of locations $P$ for computing initial loss and interior loss.
  \begin{table}[H]
      \centering
      \caption{Training data. }
      \begin{tabular}{c|ccc}
          \toprule
          Case & $N$ & \makecell{\# Sensors $m$ \\ (\# initial, \# source)}  & \makecell{\# locations $P$\\ (\# initial + \# interior)} \\
          \hline
          1D wave Eq. (case 1) & 5000 & 51,225 & 101+$51^2$ \\
          1D wave Eq. (case 2) & 5000 & 51,225 & 101+$51^2$\\
          1D wave Eq. (case 3) & 1000 & 51,225 & 101+$51^2$\\
          1D wave Eq. (case 4) & 2000 & 51,225 & 101+$51^2$\\
          2D wave Eq. & 2000 & 100,1000 & 21+$21^3$\\
          Burgers' Eq. (case 1) & 5000 & 51,225 & 101+$51^2$ \\
          Burgers' Eq. (case 2) & 5000 & 51,225 & 51+$51^2$\\
          KdV Eq. & 1000 & 51,225 & 51+$51^2$ \\
          Schr\"odinger Eq. & 5000 & 51,225 & 101+$51^2$ \\
          \bottomrule
      \end{tabular}
      \label{appendix_data}
  \end{table}
  
  \subsection*{A.2~~MIONet architectures}
  Table \ref{appendix_arch} shows the MIONet architectures used in all examples.
  The dimension of the input layer is determined by the data size, while the size of all network output layers remains the same. This consistency is essential as it enables the Hadamard product of each network's output to be taken to obtain the final output.
  \begin{table}[H]
      \centering
      \caption{MIONet architectures. The trunk net and branch nets are composed of fully connected networks. Branch-1 net is for the initial function $\varphi$, and Branch-2 net is for the source function $f$. $D_{in}$ and $D_{out}$ are the input dimension and output dimension, respectively. The ``width" and ``depth" refer to the size of the hidden layer.}
      \begin{tabular}{c|ccc}
          \toprule
          Case & \makecell{Trunk \\ ($D_{in}$-width*depth-$D_{out}$)} & \makecell{Branch-1  \\ ($D_{in}$-width*depth-$D_{out}$)} & \makecell{Branch-2 \\ ($D_{in}$-width*depth-$D_{out}$)} \\
          \hline
          1D wave Eq. (case 1) & 2-400*4-100 & 51-100*3-100 & 225-225*3-100 \\
          1D wave Eq. (case 2) & 2-400*3-100 & 51-100*3-100 & 225-225*3-100 \\
          1D wave Eq. (case 3) & 2-200*4-100 & 51-100*3-100 & 225-225*3-100 \\
          1D wave Eq. (case 4) & 2-100*3-100 & 51-100*3-100 & 225-225*3-100 \\
          2D wave Eq. & 3-200*4-200 & 100-200*3-200 & 1000-200*3-200\\
          Burgers' Eq. (case 1) & 2-400*3-100 & 51-100*3-100 & 225-225*3-100  \\
          Burgers' Eq. (case 2) & 2-400*3-100 & 51-100*3-100 & 225-225*3-100 \\
          KdV Eq. & 2-100*4-100 & 51-100*3-100 & 225-100*3-100 \\
          Schr\"odinger Eq. & 2-400*3-100 & 51-51*3-100 & 225-225*3-100  \\
          \bottomrule
      \end{tabular}
      \label{appendix_arch}
  \end{table}

  \subsection*{A.3~~Learning rate}
  The learning rate directly affects the effect of model training. 
  MIONet models are trained via mini-batch gradient descent using the Adam optimizer with default settings and exponential learning rate decay with a decay rate of $\gamma$ every $M$ iterations. Table \ref{appendix_lr} gives the base learning rate (lr) $\eta$, step size $M$ and decay rate $\gamma$ used in all examples.
  \begin{table}[H]
      \centering
      \caption{Learning rate.}
      \begin{tabular}{c|ccc}
          \toprule
          Case & base lr $\eta$ & step size $M$ & decay rate $\gamma$ \\
          \hline
          1D wave Eq. (case 1) & 0.001 & 500 & 0.96 \\
          1D wave Eq. (case 2) & 0.001 & 500 & 0.97 \\
          1D wave Eq. (case 3) & 0.001 & 500 & 0.96 \\
          1D wave Eq. (case 4) & 0.001 & 500 & 0.96 \\
          2D wave Eq. & 0.0005 & 1000 & 0.96 \\
          Burgers' Eq. (case 1) & 0.0005 & 1000 & 0.95\\
          Burgers' Eq. (case 2) & 0.0005 & 1000 & 0.96\\
          KdV Eq. & 0.001 & 1000 & 0.96 \\
          Schr\"odinger Eq. & 0.001 & 500 & 0.96 \\
          \bottomrule
      \end{tabular}
      \label{appendix_lr}
  \end{table}
  
  \subsection*{A.4~~Parameters of generating training data}
  The parameters in the analytical expression are randomly generated, with each parameter being selected from a normal distribution. The mean and variance of the normal distribution for each parameter in each calculation example are provided in Table \ref{par_gen_data_random}.

  To select the data for training, we need to calculate the errors between the generated initial values, source terms and those of the target PDE, respectively. Therefore, we need to select the domain and sampling points for calculating the errors. Here the time domain is always consistent with the computational domain of the target PDE, while the spatial domain can be selected larger than the computational domain of interest to better capture the global information of the target functions. Table \ref{par_gen_data} provides the spatial domain, number of grid points (uniformly), and tolerance $\epsilon$ for each example.
  
  \begin{table}[H]
      \centering
      \caption{Random rule.}
      \begin{tabular}{c|llll}
          \toprule
          Case & \multicolumn{4}{c}{Parameters} \\
          \hline
          \multirow{2}{*}{1D wave Eq. (case 1)} & $A\sim\N(1,1)$ & $k\sim\N(0,1)$ & $a\sim\N(0,1)$ & $w\sim\N(-2,1)$ \\
          & $b\sim\N(0,1)$ & & & \\
          \hdashline
          1D wave Eq. (case 2) & $A\sim\N(0.2,1)$ & $ a\sim\N(1,1)$ & $ b\sim\N(0,1)$ & $ \sigma\sim\N(0.2,0.5)$ \\
          \hdashline
          \multirow{2}{*}{1D wave Eq. (case 3)} & $A\sim\N(1,0.5)$ & $ a\sim\N(0.05,0.1)$ & $ c\sim\N(-1,0.5)$ & $ k\sim\N(0,0.5)$ \\
          &  $ b\sim\N(0,0.5)$ & $ d\sim\N(1,1)$ &  & \\
          \hdashline
          \multirow{2}{*}{1D wave Eq. (case 4)} & $A\sim\N(0,1)$ & $ k_1\sim\N(0,1)$ & $ w_1\sim\N(0,1)$ & $ b_1\sim\N(1,1)$ \\
          &  $ k_2\sim\N(1,1)$ & $ w_2\sim\N(1,1)$ & $ b_2\sim\N(-1,1)$ & \\
          \hdashline
          \multirow{2}{*}{2D wave Eq.} & $A\sim\N(1,1)$ & $ a_1\sim\N(2,1)$ & $ a_2\sim\N(1,1)$ & $ \sigma\sim\N(1,1)$ \\
          &  $ k_1\sim\N(0,1)$ & $ k_2\sim\N(0,1)$ & $ w\sim\N(3,1)$  & \\
          \hdashline
          \multirow{2}{*}{Burgers' Eq. (case 1)} & $A\sim\N(0.2,1)$ & $ k_1\sim\N(0.8,1)$ & $ k_2\sim\N(1,2)$ & $ w\sim\N(1,1)$ \\
          & $ b\sim\N(0,1)$ & & &\\
          \hdashline
          \multirow{2}{*}{Burgers' Eq. (case 2)} & $A\sim\N(0.2,1)$ & $ k_1\sim\N(0.8,1)$ & $ k_2\sim\N(1,2)$ & $ w\sim\N(1,1)$ \\
          & $ b\sim\N(0,1)$ & & & \\
          \hdashline
          \multirow{2}{*}{KdV Eq.} & $A\sim\N(1,1) $ & $ a\sim\N(-1,1)$ & $ c_1\sim\N(0,1)$ & $ c_2\sim\N(-1,1)$ \\
          & $ k\sim\N(0,1)$ & $ w\sim\N(0,1)$  &  &  \\
          \hdashline
          Schr\"odinger Eq. & $A\sim\N(0.5,0.5)$ & $ \zeta\sim\N(0.3,0.5)$ & $ k\sim\N(1,0.5)$ & \\
          \bottomrule
      \end{tabular}
      \label{par_gen_data_random}
  \end{table}

  \begin{table}[H]
    \centering
    \caption{Parameters of generating training data.}
    \begin{tabular}{c|ccc}
        \toprule
        Case & spatial domain & grid (initial, source) & tolerance $\epsilon$ \\
        \hline
        1D wave Eq. (case 1) & $[-3,3]$ & 51, $101\times101$ & $(\varphi_0,\varphi_1,f):(0.5,1,1)$ \\
        1D wave Eq. (case 2) & $[-2,2]$ & 51, $101\times101$ & $(\varphi_0,\varphi_1,f):(0,0,0.6)$ \\
        1D wave Eq. (case 3) & $[-1,1]$ & 101, $101\times101$ & $(\varphi_0,\varphi_1,f):(0.8,1,1)$ \\
        1D wave Eq. (case 4) & $[-1,1]$ & 51, $101\times101$ & $(\varphi_0,\varphi_1,f):(1,1,1)$ \\
        2D wave Eq. & $[-5,5]^2$ & $41\times41$, $41\times41\times41$ & $(\varphi_0,\varphi_1,f):(0.8,0.8,0.8)$ \\
        Burgers' Eq. (case 1) & $[-2,2]$ & 51, $51\times41$ &$(\varphi,f): (0,1)$ \\
        Burgers' Eq. (case 2) & $[-4,4]$ & 51, $51\times51$ & $(\varphi,f): (0,1)$\\
        KdV Eq. & $[-5,5]$ & 51, $51\times51$ & $(\varphi,f): (0.8,0.8)$ \\
        Schr\"odinger Eq. & $[-5,5]$ & 51, $101\times101$ & $(\Re(\varphi),\Im(\varphi), \Re(f), \Im(f)): (0.4,0.4,0,0)$\\
        \bottomrule
    \end{tabular}
    \label{par_gen_data}
\end{table}

\section*{B~~Impact of generating data domain}
  We investigate the influence of the spatial domain size on the results by varying the size of the domain where the data is generated. We specifically focus on case 1 of the 1D wave equation \eqref{wave} and case 1 of Burgers' equation \eqref{Burgers} as examples for our analysis. 
  
\begin{figure}[H]
      \centering
          \includegraphics[width=0.66\linewidth]{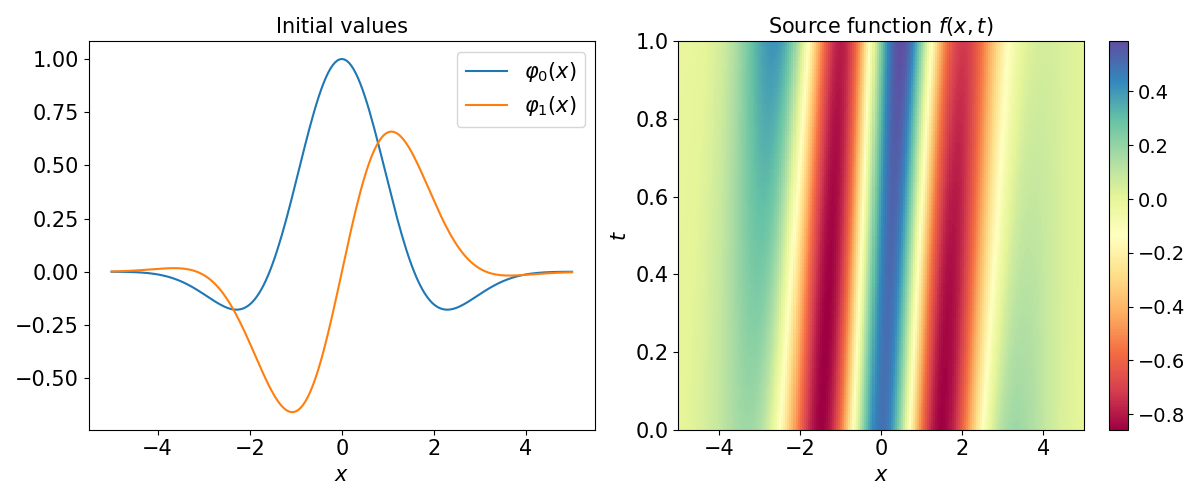}
          \includegraphics[width=0.33\linewidth]{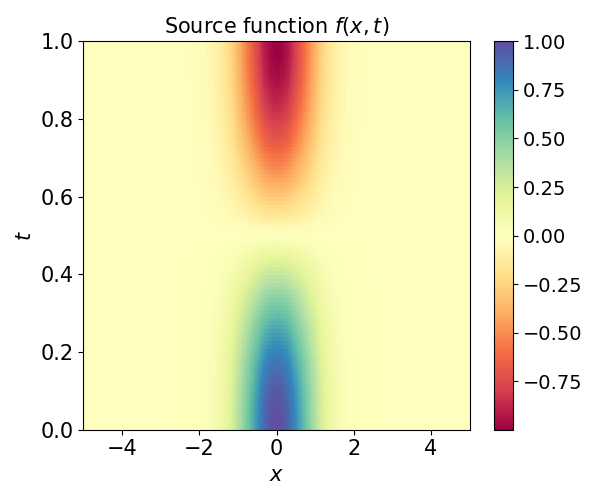}
      \caption{Left: Initial values of the target wave equation. Middle: Source function of the target wave equation. Right: Source function of the target Burgers' equation.}
      \label{test_fun}
  \end{figure}

\begin{figure}[H]
      \centering
          \includegraphics[width=0.24\linewidth]{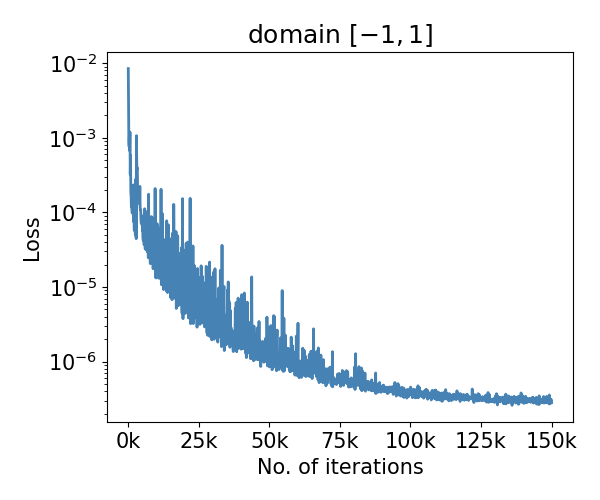}
          \includegraphics[width=0.24\linewidth]{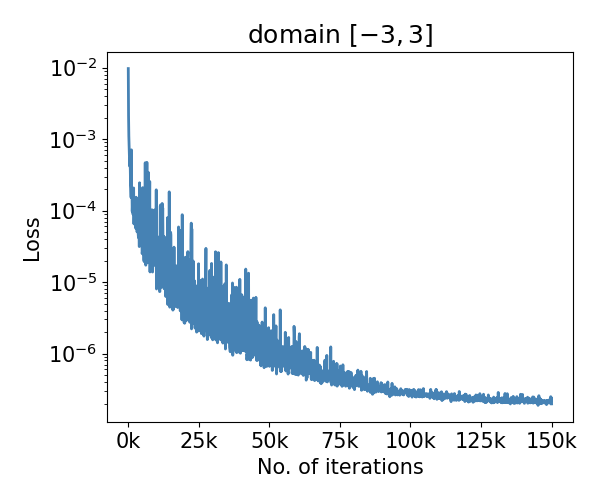}
          \includegraphics[width=0.24\linewidth]{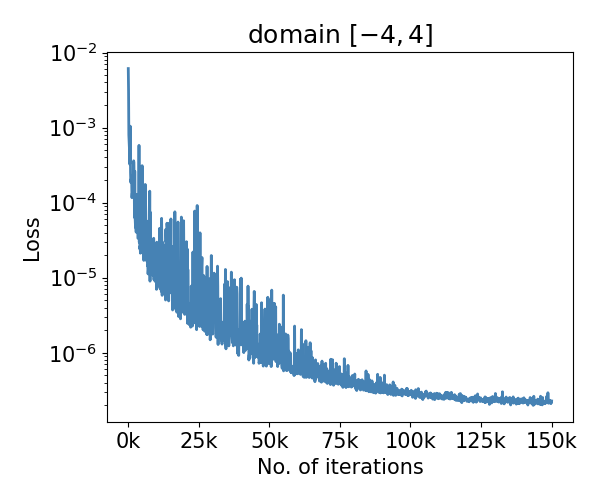}
          \includegraphics[width=0.24\linewidth]{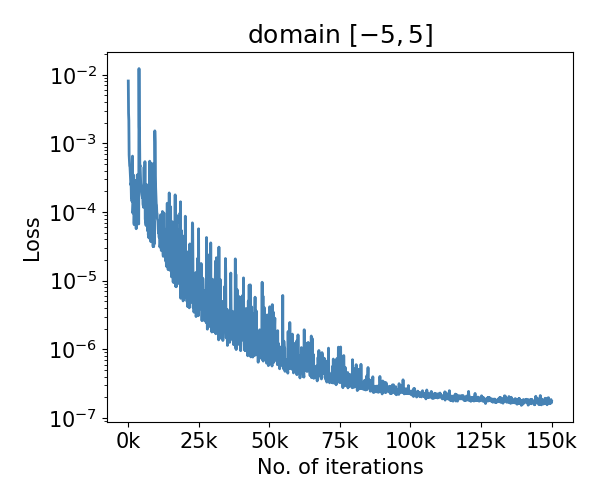} \\
          \includegraphics[width=0.24\linewidth]{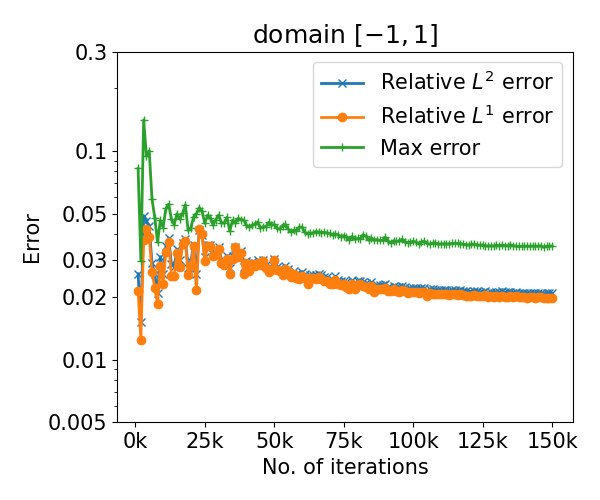}
          \includegraphics[width=0.24\linewidth]{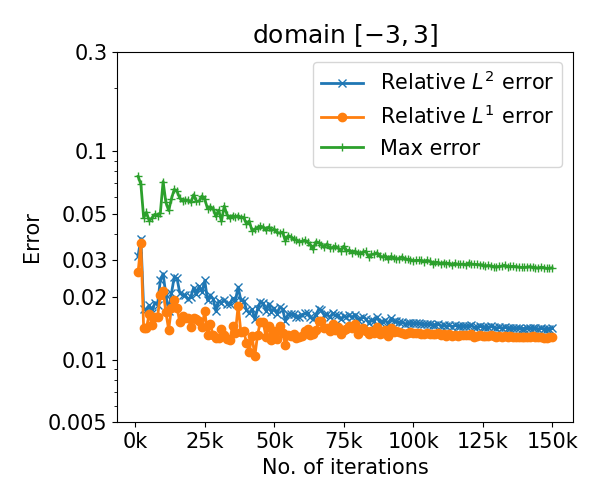}
          \includegraphics[width=0.24\linewidth]{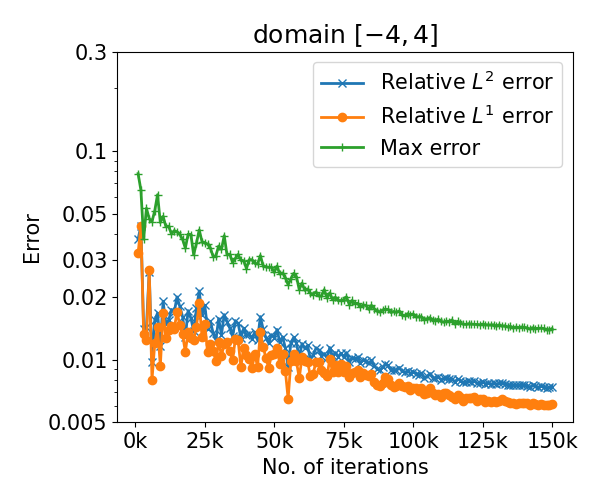}
          \includegraphics[width=0.24\linewidth]{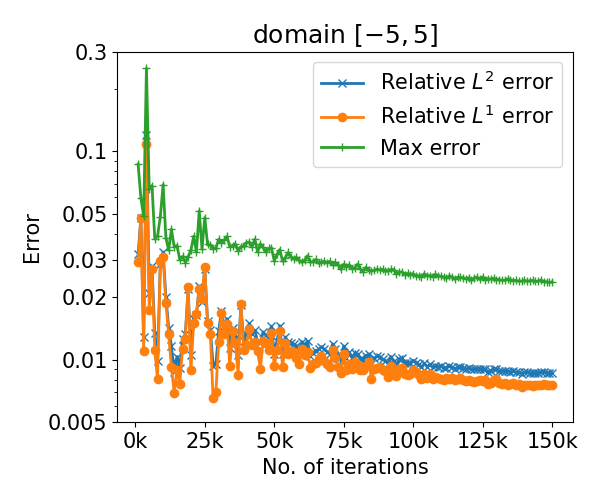}
      \caption{\textit{Solving a 1D wave equation using data generated from different domains:} Training loss and test errors for 150,000 iterations.}
      \label{wave_domain}
  \end{figure}

  To make the compact support domain of the target function larger, we introduce a slight modification to the target function \eqref{wave_case1_target} of the wave equation. The modified target functions are given as 
  \beqs
  \begin{aligned}
  &\varphi_0(x) = \exp\left(-x^2/4\right)\cos(x),\\
  &\varphi_{1}(x) = \exp\left(-x^2/4\right)\sin(x),\\
  &f(x,t) = \exp(-x/4^2)\left((-x^2/4+1/2)\cos(t-x)+x\sin(t-x)\right),
  \end{aligned}
  \eeqs
  with the exact solution $u(x,t) = \exp\left(-x^2/4\right)\cos(t-x)$. The initial values and source function are depicted in the left and right planes of Figure \ref{test_fun}, respectively. The domain of interest is defined as $[-1,1]$. We conduct tests by generating data domains of $[-1,1], [-2,2], [-3,3], [-4,4]$. Within each domain, we uniformly sample points while maintaining the same grid density, i.e., using a fixed spatial grid step size of 0.02.
  The training process and prediction errors are shown in Figure \ref{wave_domain}. The final errors are summarized in Table \ref{domain}. It can be observed that the best performance is achieved when the data domain is set as $[-4,4]$. 
  
\begin{figure}[H]
      \centering
          \includegraphics[width=0.24\linewidth]{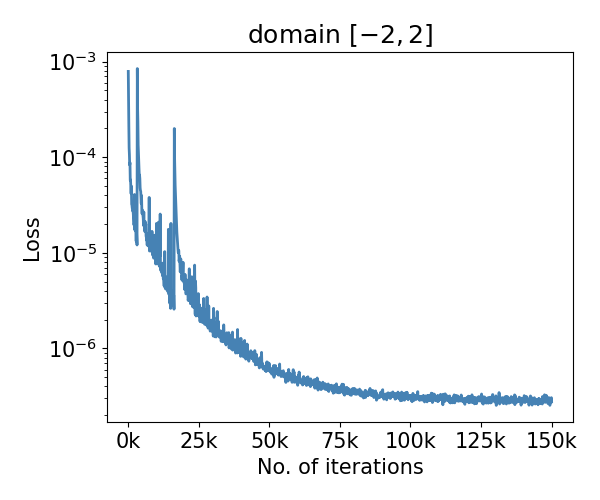}
          \includegraphics[width=0.24\linewidth]{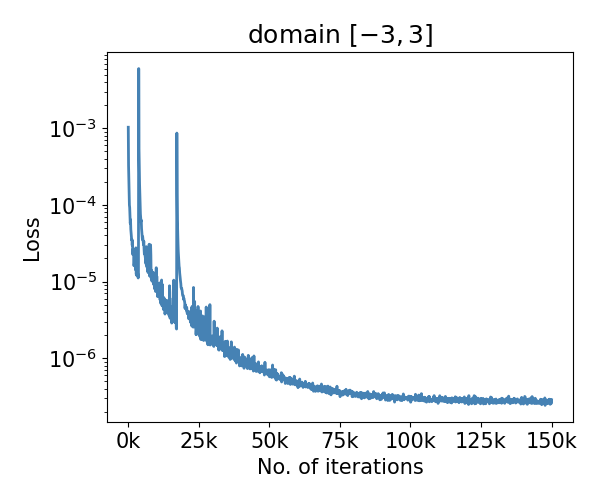}
          \includegraphics[width=0.24\linewidth]{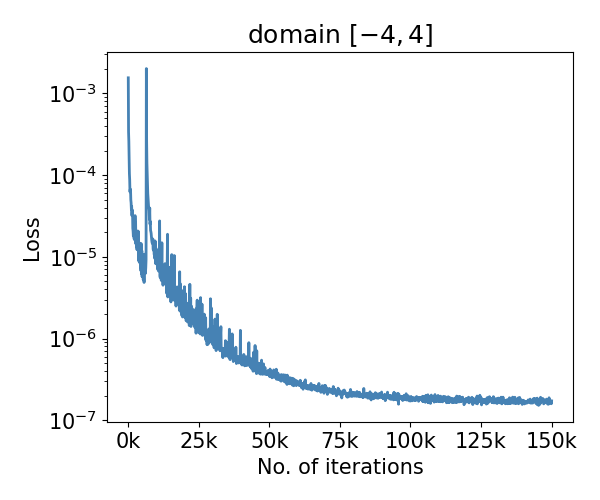}
          \includegraphics[width=0.24\linewidth]{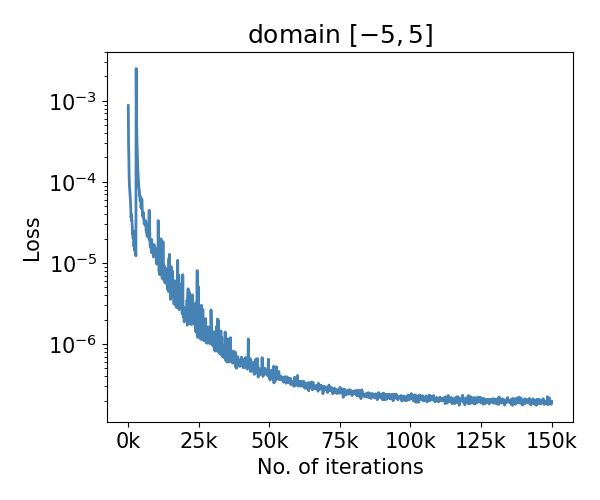} \\
          \includegraphics[width=0.24\linewidth]{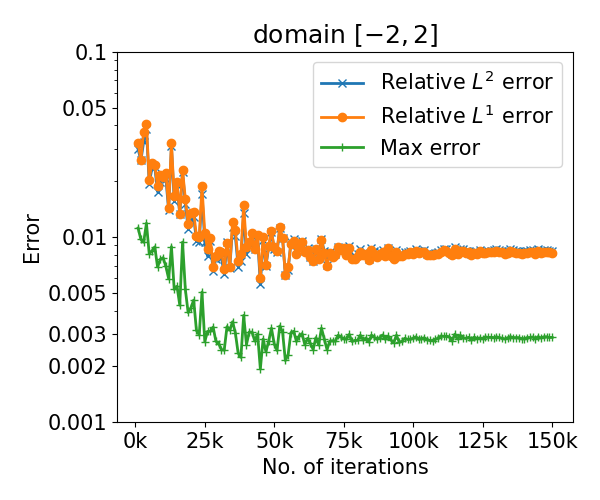}
          \includegraphics[width=0.24\linewidth]{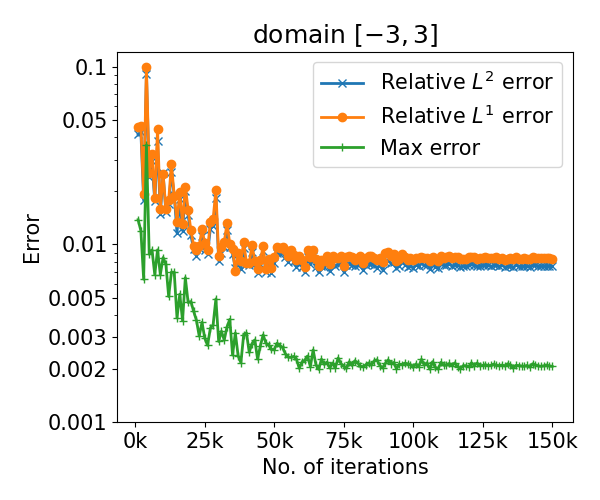}
          \includegraphics[width=0.24\linewidth]{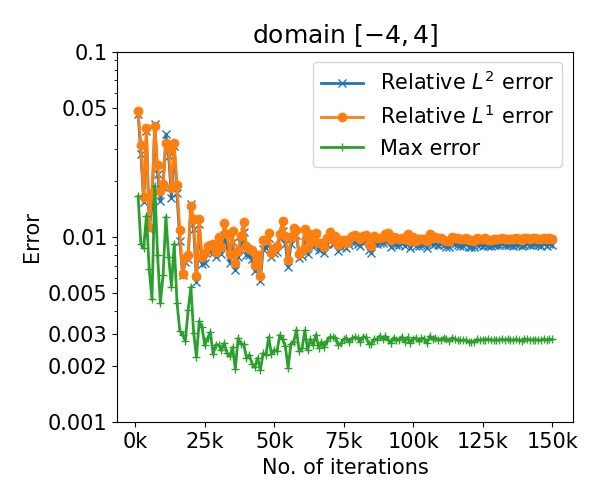}
          \includegraphics[width=0.24\linewidth]{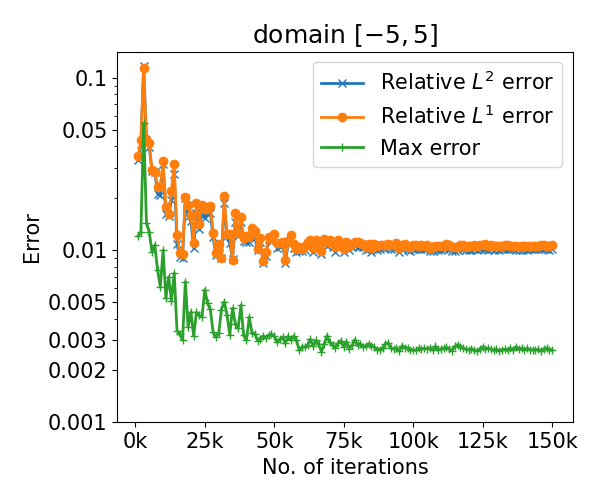}
      \caption{\textit{Solving a 1D Burgers' equation using data generated from different domains:} Training loss and test errors for 150,000 iterations.}
      \label{burgers_domain}
  \end{figure}

  For the case of Burgers' equation, we adopt the same settings as described in section \ref{sec_Burgers1}. The initial value is zero, and the source function is depicted in the right plane of Figure \ref{test_fun}. We generate data on domains $[-2,2], [-3,3], [-4,4], [-5,5]$ to test the performance. 
  We uniformly collect points within the domain and maintain the same grid density. That is, we fix the spatial grid step size as 0.02.
  The training process and prediction errors are shown in Figure \ref{burgers_domain}. The final errors are shown in Table \ref{domain}, which indicates that the results are optimal when the computational domain is chosen as $[-3, 3]$.

  Combining these two examples, we can tentatively conclude that a larger computational domain does not necessarily guarantee better results. However, the selected computational domain should not be too small either. It is preferable to choose a domain that includes the compact support set of target functions,  meaning that target functions decay toward zero within that domain.
   
  \begin{table}[H]
      \centering
      \caption{Impact of generating data domain.} 
      \begin{tabular}{c|ccc|c|ccc}
          \toprule
           \multicolumn{4}{c|}{1D wave Eq. (case 1)} & \multicolumn{4}{c}{Burgers' Eq.} \\
          \hline
          domain & $L^2$-error & $L^1$-error & max-error & domain & $L^2$-error & $L^1$-error & max-error \\
          \midrule
          $[-1,1]$ & 2.08e-02 & 1.98e-02 & 3.49e-02 & $[-2,2]$ & 8.42e-03 & 8.24e-03 & 2.86e-03 \\
          $[-3,3]$ & 1.40e-02 & 1.28e-02 & 2.76e-02 & $[-3,3]$ & \textbf{7.56e-03} & \textbf{8.33e-03} & \textbf{2.08e-03} \\
          $[-4,4]$ & \textbf{7.40e-03} & \textbf{6.08e-03} & \textbf{1.40e-02} & $[-4,4]$ & 9.05e-03 & 9.80e-03 & 2.78e-03 \\
          $[-5,5]$ & 8.64e-03 & 7.54e-03 & 2.37e-02 & $[-5,5]$ & 1.06e-02 & 1.05e-02 & 2.64e-03 \\
          \bottomrule
        \end{tabular}
        \label{domain}
  \end{table}

  \section*{C~~Impact of the number of input functions $N$}
  To investigate the impact of the number of input functions $N$ on the final prediction, we choose different $N$ values to train the model for the example of Burgers' equation and KdV equation in section \ref{sec_Burgers1} and \ref{sec_KdV}, respectively. The values of hyperparameters remain unchanged. 
  
  The experimental results suggest that the prediction accuracy is not significantly influenced by $N$.
  However, as we only tested it on one example, the prediction accuracy is more sensitive to the model parameters. Different model initialization or learning rates can lead to substantial fluctuations in the final prediction.
  
  \begin{table}[H]
  \caption{Impact of $N$.} 
  \centering
  \begin{tabular}{c|ccc|ccc}
      \toprule
       &\multicolumn{3}{c|}{Burgers' Eq.} & \multicolumn{3}{c}{KdV Eq.} \\
      \cline{2-7}
      $N$ & $L^2$-error & $L^1$-error & max-error & $L^2$-error & $L^1$-error & max-error \\
      \midrule
      500 &7.60e-03 &8.46e-03 &2.02e-03 & 1.90e-03& 1.82e-03 & 3.63e-03\\
      1000 &6.28e-03 &6.64e-03 &2.39e-03 & 1.34e-03 &1.30e-03 &3.04e-03 \\
      5000 &7.17e-03 &7.74e-03 &2.39e-03 &0.84e-03 &0.83e-03 &1.38e-03 \\
      10000 &\textbf{4.88e-03} &\textbf{5.35e-03} &\textbf{1.75e-03} & \textbf{0.41e-03} & \textbf{0.37e-03} & \textbf{1.14e-03} \\
      \bottomrule
    \end{tabular}
  \end{table}
  
  \begin{figure}[H]
      \centering
          \includegraphics[width=0.24\linewidth]{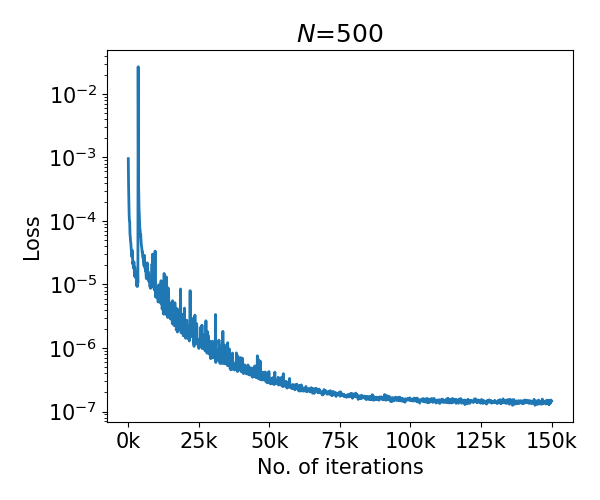}
          \includegraphics[width=0.24\linewidth]{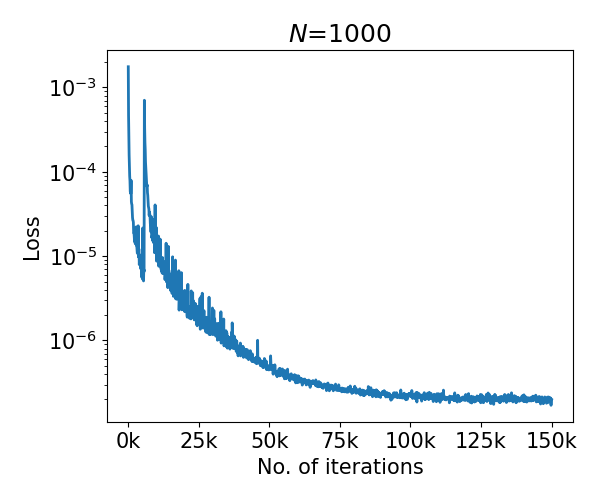}
          \includegraphics[width=0.24\linewidth]{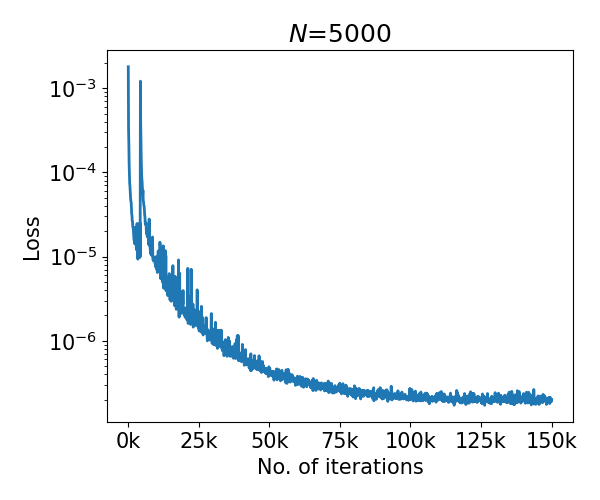}
          \includegraphics[width=0.24\linewidth]{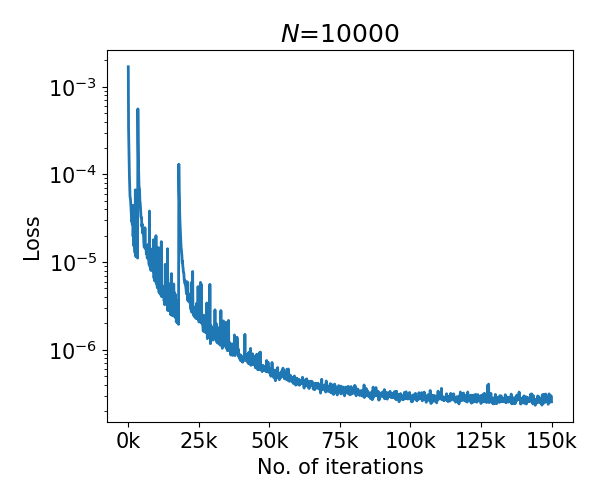} \\
          \includegraphics[width=0.24\linewidth]{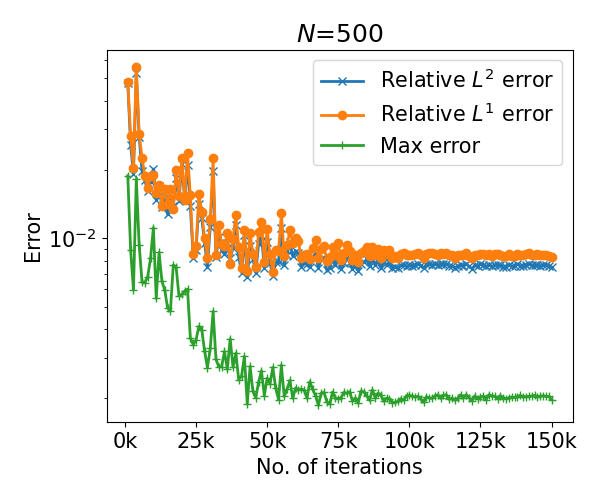}
          \includegraphics[width=0.24\linewidth]{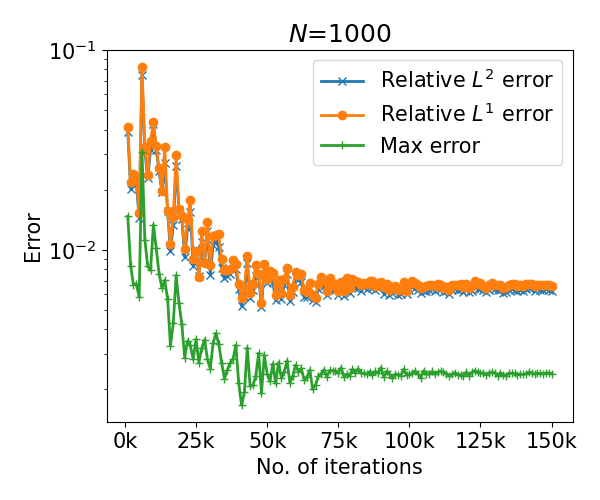}
          \includegraphics[width=0.24\linewidth]{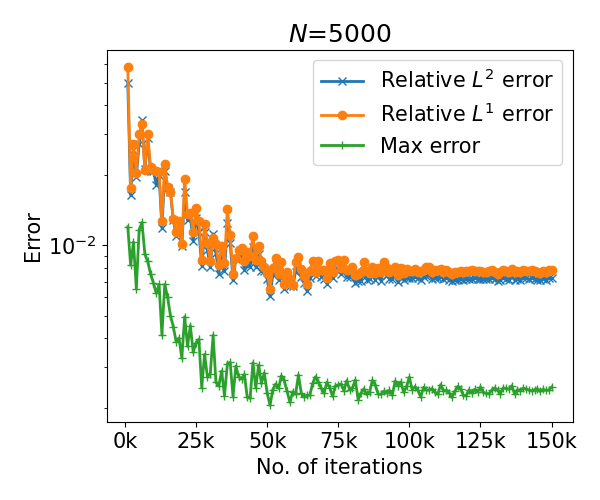}
          \includegraphics[width=0.24\linewidth]{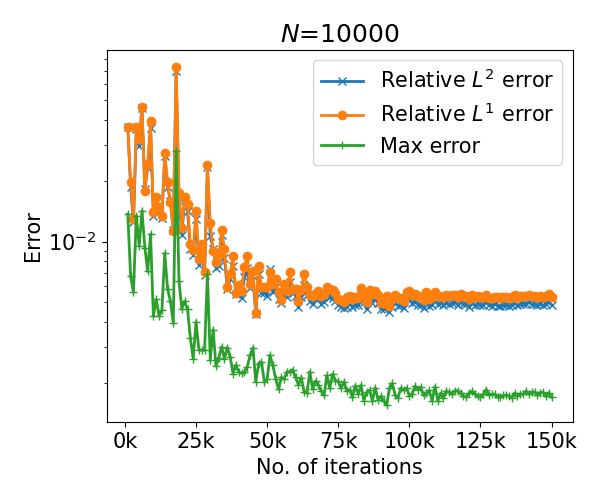}
      \caption{\textit{Burgers' equation:} Training loss and test errors for 150,000 iterations under $N=500,1000,5000,10000$.}
  \end{figure}
  
  \begin{figure}[H]
  \centering
      \includegraphics[width=0.24\linewidth]{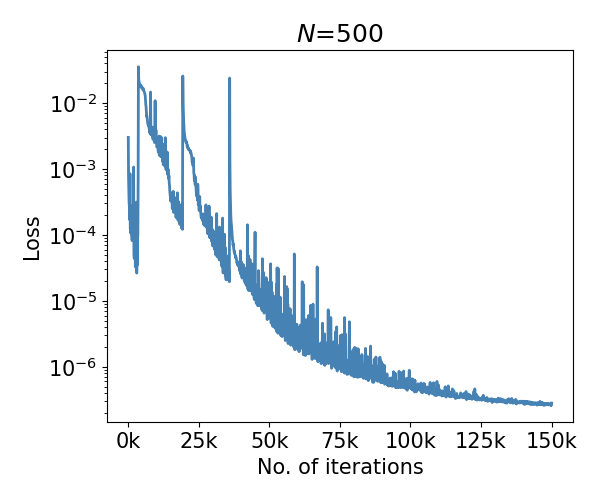}
      \includegraphics[width=0.24\linewidth]{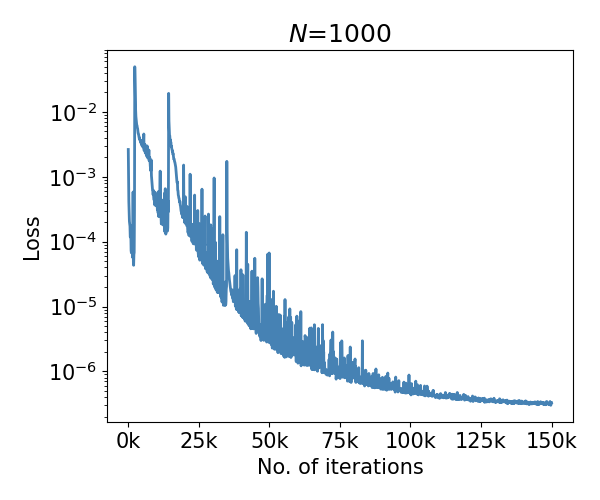}
      \includegraphics[width=0.24\linewidth]{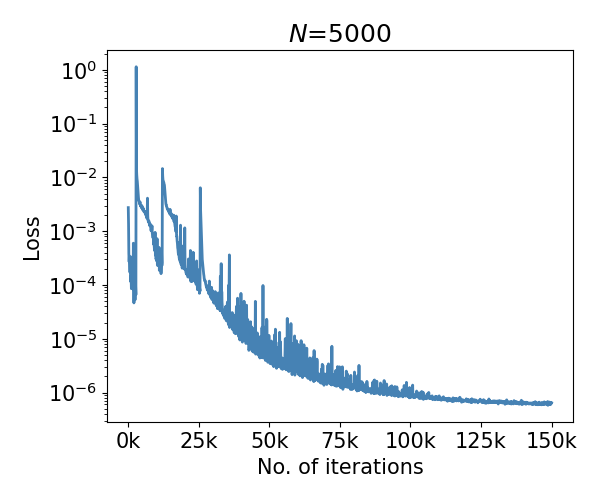}
      \includegraphics[width=0.24\linewidth]{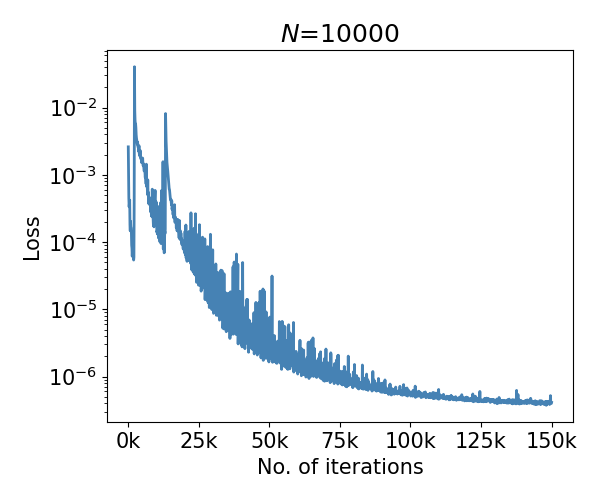} \\
      \includegraphics[width=0.24\linewidth]{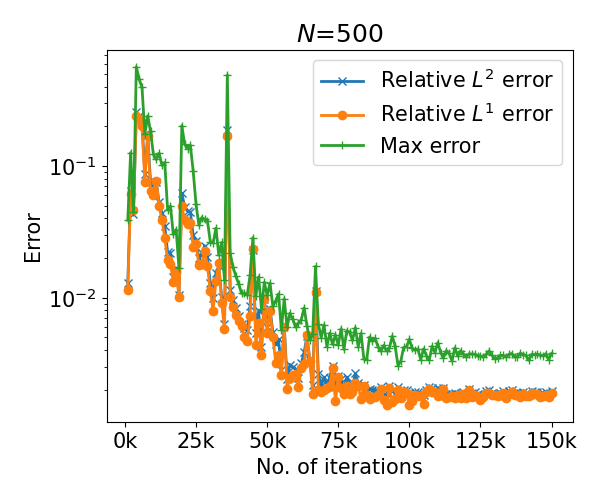}
      \includegraphics[width=0.24\linewidth]{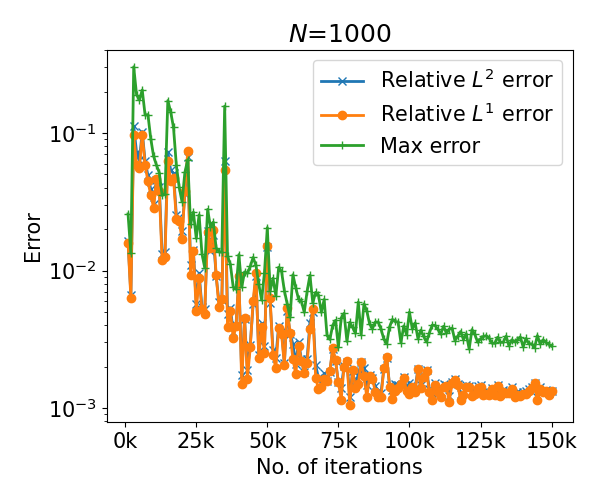}
      \includegraphics[width=0.24\linewidth]{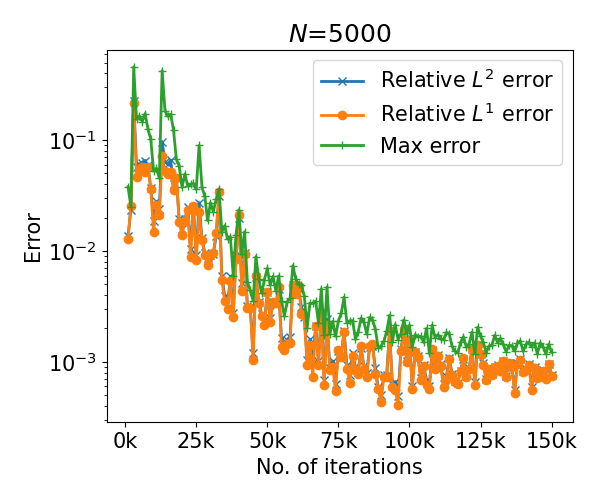}
      \includegraphics[width=0.24\linewidth]{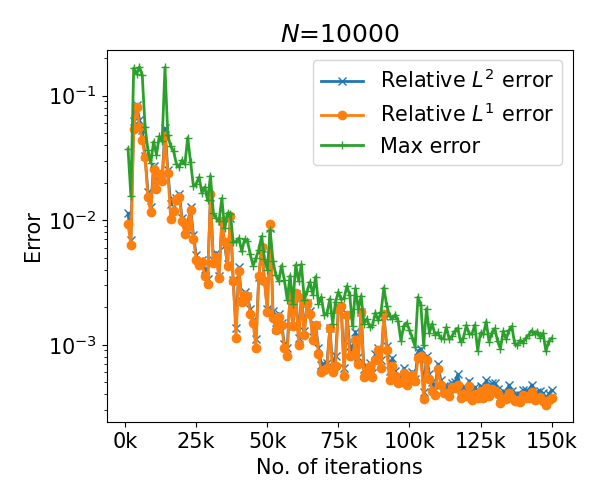}
  \caption{\textit{KdV equation:} Training loss and test errors for 150,000 iterations under $N=500,1000,5000,10000$.}
  \end{figure}

\section*{D~~Details of the experiment shown in Figure \ref{PINN}}
\begin{figure}[H]
  \centering
  \includegraphics[width=\linewidth]{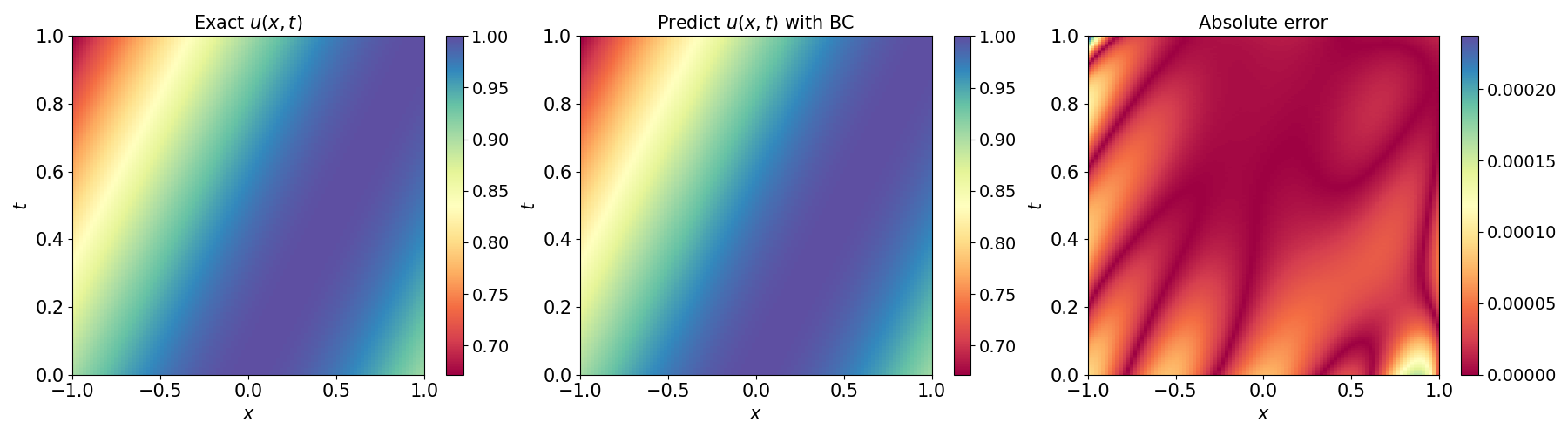}\\
  \includegraphics[width=\linewidth]{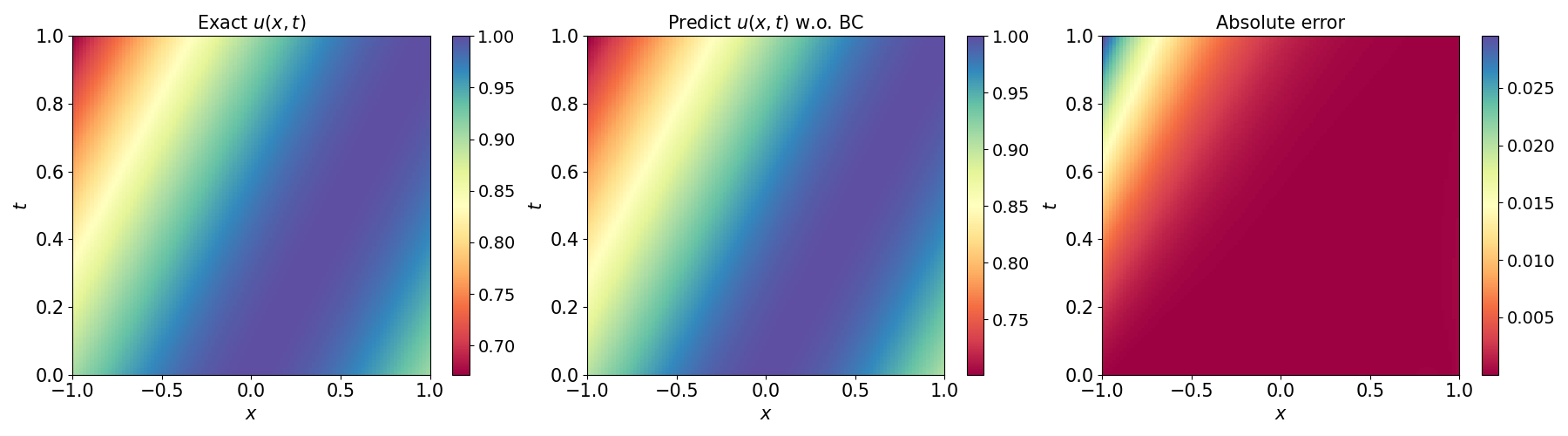}
  \caption{\textit{Solving a 1D Burgers' equation using PINNs:} Top: Training with the boundary conditions. Bottom: Training without the boundary conditions.}
\label{PINN_Burgers}
\end{figure}

To show the limitation of the PINNs in solving PDEs on unbounded domains, we use the PINNs to solve a 1D Burgers' equation
\begin{eqnarray}
  \begin{aligned}
  & u_{t}+uu_{x}-0.1 u_{xx} = f(x,t), && x\in \mathbb{R},~ t\in(0,1], \\
  & u(x, 0) = \varphi(x), && x\in \mathbb{R},\\
  \end{aligned}
\end{eqnarray}
Consider the scenario with an exact solution $u(x,t) = \exp(-0.1(x-t)^2)$. The spatial domain of interest is $\Omega=[-1,1]$. We sampled points solely at the initial time step and within the domain to train the model. As a comparative analysis,  additional experiments were conducted by sampling points not only at the initial time step and within the domain but also at the boundaries. All other settings for the two sets of experiments were consistent. We utilized fully connected neural networks with 8 hidden layers, each comprising 20 neurons, and employed the tanh activation function. The number of sampling points at the initial time step and on each boundary is set to 101, while 1000 points are sampled within the interior of the domain. The computational results across the entire domain are depicted in the Figure \ref{PINN_Burgers}. It is evident that when trained without using the boundary conditions, a notable overall error is observed.

\end{document}